\documentclass[12pt,twoside]{article}
\usepackage{amsmath,amssymb,mathrsfs,graphicx}
\usepackage[small,sc]{caption}
\usepackage{bbm}

\newtheorem{theorem}{Theorem}[section]
\newtheorem{lemma}[theorem]{Lemma}
\newtheorem{proposition}[theorem]{Proposition}
\newtheorem{corollary}[theorem]{Corollary}

\newtheorem{rmrk}[theorem]{Remark}

\DeclareMathAlphabet{\mathbfit}{OML}{cmm}{b}{it}

\makeatletter
\@addtoreset{equation}{section}
\makeatother

\newcommand{\fig}[3] {
\medskip\smallskip
\begin{figure}[htb]
  \centering
  \includegraphics[width=#2]{#1.pdf}
  \begin{minipage}[t]{0.80\linewidth}
    \caption{#3}
    \protect\label{#1}
  \end{minipage}
\end{figure}
\medskip
}

\newcommand{\figsimul}[2] {
\medskip\smallskip
\begin{figure}
  \centering
  \includegraphics[width=10cm]{s_x0_65_#1.pdf}
  \begin{minipage}[t]{0.80\linewidth}
    \caption{#2}
    \protect\label{#1}
  \end{minipage}
\end{figure}
\medskip
}

\newenvironment{remark}
{\begin{rmrk} \em}
{\end{rmrk}}

\newcommand{\fn} {function}
\newcommand{\me} {measure}

\newcommand{\erg} {ergodic}
\newcommand{\sy} {system}

\newcommand{\pr} {probability}
\newcommand{\dsy} {dynamical system}
\renewcommand{\o} {orbit}
\newcommand{\R} {\mathbb{R}}
\newcommand{\C} {\mathbb{C}}

\newcommand{\Z} {\mathbb{Z}}
\newcommand{\N} {\mathbb{N}}
\newcommand{\qed} {\hfill {\small Q.E.D.} \par\medskip}
\newcommand{\skippar} {\par\medskip}
\newcommand{\ds} {\displaystyle}
\newcommand{\proof} {\noindent \textsc{Proof.} }
\newcommand{\proofof}[1] {\noindent \textsc{Proof of {#1}.} }
\newcommand{\article}[3] {\textsc{{#1}}, {\itshape {#2}}, {{#3}}.}
\newcommand{\book}[3] {\textsc{{#1}}, {\itshape {#2}}, {{#3}}.}
\newcommand{\vol} {\textbf}
\newcommand{\eps} {\varepsilon}

\newcommand{\rset}[2] {\left\{ #1 \: \left| \: #2 \right. \! \right\} }
\newcommand{\lset}[2] {\left\{ \left. \! #1 \: \right| \: #2 \right\} }

\newcommand{\into} {\longrightarrow}

\renewcommand{\emptyset} {\varnothing}

\newcommand{\ob} {observable}
\newcommand{\ps} {X}
\newcommand{\ba} {\mathcal{A}}   % Birkhoff average
\newcommand{\bs} {\mathcal{S}}   % Birkhoff sum
\newcommand{\sa} {\mathscr{A}}   % generic sigma-algebra
\newcommand{\leb} {m}   % Lebesgue measure
\newcommand{\tailk} {X_K}   % tail of partition
\newcommand{\tow} {\mathcal{T}}   % tower map (for appendix)

\def\var{\varepsilon}

\begin{document}

\title{\textbf{Pointwise convergence of Birkhoff averages for
  global observables}}

\author{
\scshape
Marco Lenci\,\thanks{
Dipartimento di Matematica, Universit\`a di Bologna,
Piazza di Porta San Donato 5, 40126 Bologna, Italy.
E-mail: \texttt{marco.lenci@unibo.it}.}
\thanks{
Istituto Nazionale di Fisica Nucleare,
Sezione di Bologna, Via Irnerio 46,
40126 Bologna, Italy.}
\,and
Sara Munday\,\thanks{
Dipartimento di Matematica, Universit\`a di Pisa,
Largo Bruno Pontecorvo 5, 56127 Pisa, Italy.
E-mail: \texttt{sara.munday@dm.unipi.it}.}
}

\date{Final version for \emph{Chaos} \\[6pt]
August 2018}

\maketitle

\begin{abstract}
  It is well-known that a strict analogue of the Birkhoff Ergodic
  Theorem in infinite ergodic theory is trivial; it states that for any
  infinite-measure-preserving ergodic system the Birkhoff average of
  every integrable function is almost everywhere zero. Nor does
  a different rescaling of the Birkhoff sum that leads to a non-degenerate
  pointwise limit exist. In this paper we give a version of Birkhoff's theorem
  for conservative, ergodic, infinite-measure-preserving dynamical
  systems where instead of integrable functions we use certain elements
  of $L^\infty$, which we generically call global observables. Our main
  theorem applies to general systems but requires an hypothesis
  of ``approximate partial averaging'' on the observables. The idea behind
  the result, however, applies to more general situations, as we show with
  an example. Finally, by means of counterexamples and numerical
  simulations, we discuss the question of finding the optimal class of
  observables for which a Birkhoff theorem holds for
  infinite-measure-preserving systems.

\end{abstract}

\bigskip\medskip
\textbf{Birkhoff's Ergodic Theorem is a cornerstone of the theory of
dynamical systems. It states that for a dynamical system endowed with a
finite invariant measure, the time, or Birkhoff, average of an integrable
function exists almost everywhere. For an ergodic system, this is
equivalent to the Strong Law of Large Numbers for the evolution of any
integrable function. When the invariant measure is infinite, which
is the case, for example, for most unbounded or extended Hamiltonian
systems, Birkhoff's theorem is no longer significant, in the sense that at
least for ergodic systems, the Birkhoff average of any integrable observable
is almost everywhere zero. However, for a dynamical system preserving an
infinite measure, the integrable functions are not the only observables of
interest. For example, for an extended Hamiltonian system, the kinetic
and potential energies, and many other ``delocalized'' observables are
not integrable. In this article we make steps towards a formulation of
Birkhoff's Ergodic Theorem for global observables in
infinite-measure-preserving systems. A global observable is, in essence,
a bounded function that is significantly different from zero throughout the
space.
}

\section{Introduction}
\label{sec-intro}

Birkhoff's Ergodic Theorem is one of the cornerstones of probability theory
and the theory of \dsy s. It states that if $T$ is a measure-preserving
transformation of
a \pr\ space $(\ps, \mu)$ and $f \in L^1(\ps, \mu)$, the
\emph{Birkhoff average}
\begin{equation} \label{ba-intro}
  \ba f(x) := \lim_{n \to \infty} \, \frac1n \sum_{k=0}^{n-1} f\circ T^k(x)
\end{equation}
exists for $\mu$-a.e.\ $x \in \ps$. If $T$ is also ergodic, the theorem states
in addition that
\begin{equation}
  \ba f(x) = \int_\ps f \, d\mu
\end{equation}
for a.e.\ $x$. Here and in the rest of the paper we use the most general
definition of \erg ity, which is valid for both finite and infinite \erg\ theory:
The map $T$ is said to be ergodic if every invariant set $B$ (this means
that $T^{-1}B = B$ mod $\mu$, where mod $\mu$ indicates that these
sets are equal up to a $\mu$-null set of points) has zero \me\ or full \me\
(so either $\mu(B)=0$ or $\mu( \ps \setminus B) = 0$).
Thus, for a probability-preserving \sy, \erg ity corresponds to the Strong
Law of Large Numbers for the variables $f \circ T^n$, for all $f \in L^1$.

Let us now consider the case where $(\ps, \mu)$ is an infinite \me\
space. More precisely, let us assume that $\mu$ is a $\sigma$-finite infinite
\me, which means that $\mu(X) = \infty$ and $X$ can be written as
$X = \bigcup_{j \in \N} X_j$, with each $\mu( X_j ) < \infty$. For an \erg\ $T$
the strict analogue of Birkhoff's Theorem is trivial: For all $f \in L^1$,
$\ba f(x) = 0$ almost everywhere. (This is an easy consequence of Hopf's
Ergodic
Theorem, which we recall momentarily.) One is thus led to ask what
the growth rate is for the \emph{Birkhoff sum}
\begin{equation}
  \bs_n f := \sum_{k=0}^{n-1} f \circ T^k
\end{equation}
of an integrable \fn\ $f$. Aaronson discovered that there is no universal
growth rate. More precisely, if $T$ is conservative (in other words, Poincar\'e
recurrence holds \cite[\S1.1]{a}) and \erg, given any sequence
$(a_n)_{n \in \N}$ of positive numbers, one of the following two cases
occurs \cite[Thm.\ 2.4.2]{a}:
\begin{enumerate}
\item For all $f \in L^1$ with $f>0$, $\ds \liminf_{n \to \infty} \,
  \frac{\bs_n f(x)} {a_n} = 0$, for a.e.\ $x \in \ps$.

\item There exists a strictly increasing sequence $(n_k)_{k \in \N}$ of
  the natural numbers such that for all $f \in L^1$ with
  $f>0$, $\ds \lim_{k \to \infty} \frac{\bs_{n_k} f(x)} {a_{n_k}} = \infty$,
  for a.e.\ $x \in \ps$.
\end{enumerate}

Notice in the second case that the sequence $(n_k)$ is the same
for all $f$ and for all $x$, and thus the assertion is stronger than the
statement: For every $f \in L^1$ with $f>0$, $\limsup_{n \to \infty}
\bs_n f / a_n = \infty$ almost everywhere.

The lack of a universal growth rate for $(\bs_n f(x))_{n \in \N}$ is not
due to its dependence on $f$, but on $x$. In fact, for an \erg\ $T$,
Hopf's Ergodic Theorem \cite{s,h} states that, for all $f,g \in L^1$
with $g > 0$,
\begin{equation}
  \lim_{n \to \infty} \frac{\bs_n f(x)} {\bs_n g(x)} = \frac{\int_\ps f \, d\mu}
  {\int_\ps g \, d\mu},
\end{equation}
for a.e.\ $x$. Thus, if we choose a function $g \in L^1$ with $g > 0$ and
$\int g \, d\mu = 1$, and set $a_n(x) := \bs_n g(x)$, we indeed have
that $\bs_n f(x) / a_n(x) \to \int f d\mu$ almost everywhere, for all
$f \in L^1$. But the variability of $x \mapsto (a_n(x))_{n \in \N}$ is
so strong that only a zero-\me\ set of points produces the same rate.

In this paper, \fn s $f: \ps \into \C$ are supposed to represent observations
about the state of the \sy\ $x \in \ps$. Accordingly, they will be called
\emph{\ob s}. In particular, all functions $f \in L^1(\ps, \mu)$ will be called
\emph{local \ob s}. The name is due to the fact that they are well approximated
by \fn s with a finite-\me\ support within an infinite-\me\ ambient space.

The results presented above lead one to think that local \ob s are not
the right ones to average along the \o s of $T$.
Even when a scaling sequence exists such that $\bs_n f / a_n$
converges to a non-degenerate limit, cf.\ the Darling-Kac Theorem
\cite[\S3.6]{a}, this convergence is in distribution (more precisely,
\emph{strongly in distribution}) and the
limit is a non-constant random variable. In any case, the limit cannot
reasonably be called the average of $f$ along the \o s of $T$. A
natural concept of average presupposes that if, for example, $f \equiv c$,
its average is $c$. In other words, we are interested in \emph{bona fide}
Birkhoff averages, as in (\ref{ba-intro}).

So we need to change the class of \ob s. The simplest class beyond
$L^1(\ps, \mu)$ that one might think to consider is $L^\infty(\ps, \mu)$,
which does include the constant \fn s. However, the whole of $L^\infty$ is,
vaguely speaking, ``too big'' for us to expect constant Birkhoff averages
for all of its elements. An interesting class of counterexamples is given
by the indicator \fn s of infinite-\me\ sets with the property that \o s spend
long stretches of time there before leaving; for example, neighborhoods
of strongly neutral indifferent fixed points. The Birkhoff averages of these
\ob s converge, strongly in distribution, to non-constant random variables.
A classical example of this phenomenon is the arcsine law for the Boole
transformation \cite{t02}. We will return to this example, along with others,
in Section \ref{sec-counter}.

In this paper we call \emph{global \ob s} all essentially bounded \fn s
for which, in principle, a Birkhoff Theorem could hold. Of course,
depending on the \sy\ at hand, the Birkhoff Theorem will hold as well
for many non-integrable, non-essentially bounded \ob s. Nonetheless, here
we limit ourselves to subspaces of $L^\infty$, for two reasons. 
Firstly, as already discussed, 
$L^\infty$ already contains ``too many'' observables. Secondly, we
want to follow the approach of Lenci on the question of mixing for
infinite-\me-preserving \dsy s \cite{limix, lpmu}, whereby global \ob s are
taken from subspaces of $L^\infty$. (This is an important assumption
there because the theory exploits the duality between $L^1$ and $L^\infty$.)
Another observation to make is that with the vague ``definition'' given
above, it is impossible to pre-determine the space of global \ob s. We do in
fact expect it to depend significantly on the given \sy. Nevertheless,
the common underlying concept can be expressed like this: a
global observable is a \fn\ which is supported more or less all over the
infinite-\me\ space and which \me s a quantity that is roughly homogeneous
in space. For example, if the reference space is $(\R^d, \leb)$, where $\leb$
is the Lebesgue \me, all periodic or quasi-periodic bounded \fn s are in
principle global \ob s. Another example is the case where $T : [0,1] \into [0,1]$
is an expanding map with an indifferent fixed point at 0 and preserves an
absolutely continuous \me\ that is non-integrable around 0. Then all
bounded \fn s which have a limit at 0 or oscillate in a controlled way in its
neighborhood are candidates for global \ob s.

\skippar

The main result of this paper, which we present in Section
\ref{sec:mainthm}, is an analogue of the Birkhoff Theorem for certain global
\ob s relative to a conservative, \erg, infinite-\me-preserving \dsy\
$(\ps, \mu, T)$. The hypotheses of the theorem are formulated in terms of
the partition of $\ps$ determined by the hitting times to a set $L_0$. This
partition is a very natural construction; for \sy s isomorphic to a Kakutani
tower, which includes all invertible maps \cite[\S 1.5]{a}, it corresponds to
the levels of the tower. In the appendix, we recall the definition and basic
properties of Kakutani towers. Returning to Section \ref{sec:mainthm},
we also describe several concrete examples of systems and \ob s for
which our results hold.

As our main theorem is certainly not optimal, we further discuss its core
ideas and limitations. First, in Section \ref{sec-alpha}, we give an example
of a family of dynamical systems --- which happen to be conjugates
of $\alpha$-Farey maps \cite{kms} --- and a family of global \ob s which do
not satisfy the hypotheses of the theorem, but for which we are nevertheless
able to prove that the Birkhoff average is almost everywhere constant. This
is done using the same ideas as in the proof of the theorem, but the
techniques are rather more complicated and specific to that case. Finally,
in Section \ref{sec-counter} we briefly recall the known examples mentioned
above of $L^\infty$ \fn s whose Birkhoff average does not converge
pointwise, and we construct other examples of a similar nature which are
interesting because they are representations of L\'evy walks (see
\cite{zdk, cgls, mssz} and references therein), thus highlighting the
connections between infinite \erg\ theory and anomalous stochastic processes.
In light of the vague definition given above, these \fn s cannot really be
considered counterexamples to our theorem. So we also present numerical
simulations of the Birkhoff averages for the observables and the systems
discussed in Sections \ref{sec:mainthm} and \ref{sec-alpha}.

Let us also remark here that there is a related strand of research in
which finite-measure spaces with non-$L^1$ observables are investigated.
For instance, in the 90s, first Major \cite{m} and then Buczolich \cite{b}
constructed specific examples where two different finite-measure systems
assign almost everywhere a different constant value to the limit of the
Birkhoff averages for the same non-$L^1$ observable. That is, the limit
$\ba f$ exists almost everywhere and is constant for each of the two systems,
but it is not equal to the integral of $f$ (see also the survey article
\cite{bsurv} for related references). More recently, Carney and Nicol
\cite{cn} investigate growth rates of Birkhoff sums of non-integrable
observables. Also, the effect on the strong law of large numbers of
``trimming'' the largest value(s) from the sums has been studied first
by Aaronson and Nakada \cite{an} and then by Kesseb\"ohmer and
Schindler \cite{ks}.

\bigskip\noindent
\textbf{Acknowledgments.}\ This research is part of the authors'
activity within the \linebreak DinAmicI community, see
\texttt{www.dinamici.org} and also part of M.L.'s activity within the
\emph{Gruppo Nazionale di Fisica
Matematica} (INdAM, Italy).

\section{Setup and main theorem}
\label{sec:mainthm}

For the rest of this paper, we will indicate a \dsy\ by means of a
triple $(\ps, \mu, T)$, where $\ps$ is a measurable space, $\mu$ a \me\
on it and $T: \ps \into \ps$ a measurable map. We shall always assume
that $(\ps, \mu)$ is a $\sigma$-finite \me\ space. Strictly speaking,
we should also mention the $\sigma$-algebra $\sa$ of all measurable
sets of $\ps$, but, as we only deal with one $\sigma$-algebra, we shall take it as understood.
 Unless otherwise stated, we shall always assume that
$T$ preserves $\mu$, meaning that, for all measurable $A \subseteq \ps$,
$\mu(T^{-1} A)  = \mu(A)$. The measure $\mu$ can be infinite, that is,
$\mu( \ps ) = \infty$, or finite, in which case we assume it to be normalized,
that is, $\mu( \ps ) = 1$. Here we are particularly interested in the first case.

So, given a \dsy\ $(\ps, \mu, T)$, an \ob\ $f: \ps \into \C$, a positive
integer $n$  and a point $x \in \ps$, we denote:
\begin{align}
  \bs_n f(x) &:= \sum_{k=0}^{n-1} f\circ T^k(x); \\[3pt]
  \ba_n f(x) &:= \frac{ \bs_n f(x) } n; \\[10pt]
  \ba f(x) &:= \lim_{n \to \infty} \ba_n f(x),
\end{align}
whenever the limit exists.

Our goal is to find conditions under which $\ba f(x)$ exists and is
constant almost everywhere. The easiest such condition is perhaps that
$f$ is a \emph{coboundary}, as in the next proposition, whose proof is
trivial.

\begin{proposition} \label{prop0}
  For a \dsy\ as described above, let $f = g - g \circ T^k$, with
  $g \in L^\infty(\ps, \mu)$ and $k \in \Z^+$. Then $\ba f(x) = 0$
  $\mu$-almost everywhere.
\end{proposition}

Another simple condition was already mentioned in the introduction. We
repeat it here for completeness.

\begin{proposition} \label{prop1}
  Suppose that $(\ps, \mu, T)$ is an infinite-\me-preserving \erg\ \dsy\
  and $f \in L^1(\ps, \mu)$. Then $\ba f(x) = 0$ almost everywhere.
\end{proposition}

\begin{corollary} \label{cor1}
  If the \ob\ $f$ is such that $f - f^* \in L^1$ for some $f^* \in \C$, then
  $\ba f(x) = f^*$ almost everywhere.
\end{corollary}

Corollary \ref{cor1} applies to a large number of \ob s which converge
to a constant ``at infinity''. For this phrase to make sense, a topology
and a notion of infinity must be defined on $\ps$. However, this fact is
very general and can be stated in a purely \me-theoretic fashion, as in the
following proposition.

\begin{proposition} \label{prop2}
  For an infinite-\me-preserving \erg\ \sy\ $(\ps, \mu, T)$, suppose that
  $f\in L^\infty(\ps, \mu)$ admits $f^*\in \C$ with the following property:
  For all $\eps > 0$, there exists a finite-\me\ set $A_\eps$ such that
  $|f(x) - f^*| \le \eps$ for every $x \in \ps \setminus A_\eps$. Then,
  for $\mu$-a.e.\ $x \in \ps$, $\ba f(x) =  f^*$.
\end{proposition}

\proof Without loss of generality, suppose that $f^*=0$ (since if not, we
can always consider the function $g:=f-f^*$).

For $\eps>0$, define the \ob\ $f_\eps := f \, 1_{\ps \setminus A_\eps}$,
where $1_A$ denotes the indicator \fn\ of the set $A$. By hypothesis,
$\| f_\eps \|_\infty \le \eps$, so it follows that $\| \ba_n f_\eps \|_\infty \le
\eps$ for all $n \in \N$. Consider now the \fn\ $f - f_\eps = f \, 1_{A_\eps}
\in L^1(\ps, \mu)$. By Proposition \ref{prop1} there exists a
full-\me\ set $B_\eps$ such that, for all $x \in B_\eps$,
\begin{equation}
  \lim_{n\to\infty} \left( \ba_n f(x) - \ba_n f_\var(x) \right) =0,
\end{equation}
whence
\begin{equation} \label{main50}
  \limsup_{n\to\infty} \left| \ba_n f(x) \right| \le \var.
\end{equation}
If we choose a sequence $\eps_i \to 0$, we conclude that, for every
$x \in \bigcap_i B_{\eps_i}$, $\ba f(x) = 0$.
\qed

Proposition \ref{prop2} is extremely general and, for that reason, also
rather weak, because it works with \ob s that are almost constant on the
overwhelming largest part of the space $\ps$. Our main theorem, which we
state after giving an \emph{ad hoc} construction,  is a  stronger result
that effectively uses the dynamics of $T$.

Assume that $T$ is conservative and ergodic. Given a set $L_0$ with
$0 < \mu(L_0) < \infty$, we have that
\begin{equation}
  \bigcup_{k\ge 0}T^{-k} L_0 = \ps \ \mathrm{mod} \ \mu,
\end{equation}
that is, $L_0$ is a sweep-out set. If we recursively define, for $k \ge 1$,
\begin{equation} \label{lk}
  L_k := \left( T^{-1} L_{k-1} \right) \setminus L_0
\end{equation}
we see that $\{ L_k \}_{k \in \N}$ forms a partition of $\ps$ mod $\mu$.
(We use the convention whereby $0 \in \N$.) By construction, $L_k$ is the
set of points whose \o\ intersects $L_0$ for the first time at the
$k^\mathrm{th}$ iteration. In other words, $\{ L_k \}$ is the partition which
consists of the level sets of the hitting time to $L_0$. By (\ref{lk}),
$\mu(L_k) \le \mu(L_{k-1})$. Also, by conservativity, $\mu(L_k) \to 0$
as $k \to \infty$. Finally $\sum_k \mu(L_k) = \mu(X) = \infty$.

Observe that if $(\ps, \mu, T)$ is isomorphic to a Kakutani tower
and $L_0$ corresponds, via the isomorphism, to the base of the tower,
then $L_k$ corresponds to the $k^\mathrm{th}$ level of the tower for all $k$.
In particular, the above construction has a clear interpretation in the case
of an invertible $T$. We refer to Appendix \ref{app-kakutani} for the definition of
a Kakutani tower and results linking Kakutani towers to Theorem
\ref{thm1}.

\begin{theorem} \label{thm1}
  Let $(\ps, \mu, T)$ be an infinite-\me-preserving, conservative, \erg\ \dsy,
  endowed with the partition $\{ L_k \}_{k \in \N}$, as described above. Let
  $f\in L^\infty(\ps, \mu)$ admit $f^*\in \C$ with the following
  property: $  \forall \var>0$, $\exists N, K\in \N \text{ such that }\forall x \in
  \bigcup_{k \ge K} L_k$,
  \begin{displaymath}
  \left| \ba_N f(x)-f^* \right| \le \var.
  \end{displaymath}
  Then for $\mu$-a.e.\
  $x \in \ps$, $\ba f(x) = f^*$.
\end{theorem}

\proof Once again, it is enough to prove the theorem in the case
$f^*=0$. Also, without loss of generality, we may assume that
\begin{equation} \label{main60}
  \frac{N \|f\|_\infty} K \le \var,
\end{equation}
otherwise we can always take a larger $K$ in the main hypothesis
of the theorem. Define the sets
\begin{align}
  \tailk &:= \bigcup_{k \ge K} L_k ; \label{tailk} \\[4pt]
  \tailk^c &:= \ps \setminus \tailk  \label{tailkc}
\end{align}
and the \ob\ $f_\var := f \, 1_{\tailk}$. The hypotheses on $T$ imply that
the \o\ of $\mu$-a.e.  $x\in \ps$ must enter $\tailk$ infinitely often. Let us fix
one such $x$ and split its \o\ into blocks which are subsets, alternately, of
$\tailk$ and $\tailk^c$. So, let $m_0 = m_0(x)$ denote the first time where
$T^{m_0}(x) \in \tailk$, and note that $m_0$ can be equal to zero (in which
case $x$ is already in $\tailk$). In other words, $T^{m_0}(x) \in \tailk$ and
$T^k(x) \in \tailk^c$ for all $k<m_0$. Denote $k_1 = k_1(x) \ge K$ the
unique integer such that $T^{m_0}(x) \in L_{k_1}$. Set $n_1 := k_1 - K$, so
that $T^{m_0+n_1-1}(x) \in L_{K}$. Now, let $m_1 = m_1(x)$ denote the length
of the following excursion in $\tailk^c$, so that $T^{m_0+n_1+m_1-1}(x) \in
L_0$ and the next \o\ point jumps back to the set $\tailk$, say to the set
$L_{k_2}$, for some $k_2 = k_2(x) \ge K$. Again, set $n_2 := k_2 - K$.
Continuing in this way, we construct two sequences $(m_j)_{j \ge 0}$ and
$(n_j)_{j \ge 1}$, where setting $M_j := m_0 + m_1 + \cdots + m_j$ and
$N_j := n_1 + n_2 + \cdots +n_j$, we have that
\begin{itemize}
  \item $m_0 \ge 0$ and $m_j\geq K$ for all $j\geq1$.
  \item $T^{M_j+N_j-1}(x)\in L_0$ and $T^{M_j+N_j}(x)\in L_{k_{j+1}}\subset
    \tailk$ for all $j\ge 1$.
  \item $f_\var(T^{M_j+N_{j+1}+i}(x))=0$ for all $j\ge 0$ and $0 \le i<m_{j+1}$.
\end{itemize}

Now fix $n \ge m_0+n_1+m_1$ for which there exists $j \in \N$ such that
$n = M_j + N_j + i$, for $0 \le i < n_{j+1}$. In other words, we consider all
sufficiently large $n \in \N$ which correspond to stopping the \o\ of $x$ during an
excursion in $\tailk$. We will treat the other values of $n$ later. The only parts of
the orbit that contribute in a non-zero way to the Birkhoff sum $\bs_n f_\var(x)$
are the excursions in $\tailk$, that is the blocks of lengths $n_1, \ldots, n_j$ and $i$.
Each of these $j+1$ blocks can be further decomposed into $p_i$ sub-blocks
of length $N$ and a remainder sub-block of length $0 \le r_i < N$. By hypothesis,
the Birkhoff sum corresponding to each sub-block, except for the
remainder sub-blocks, is bounded in modulus by $\var$. The contribution of each
remainder sub-block is instead bounded by $N \| f_\var \|_\infty \le
N \| f \|_\infty $. Putting all these observations together, we have that
\begin{equation} \label{main70}
\begin{split}
  \left| \ba_n f_\var(x) \right| &< \frac1 {M_j+N_j+i} \left( \var \sum_{i=1}^{j+1}
  p_i + (j+1) N \| f \|_\infty \right) \\
  &< \var \frac{ \sum_{i=1}^{j+1} p_i } { M_j+N_j+i } + \frac{ (j+1) N \|f\|_\infty}
  {jK} \\[4pt]
  &< \var + 2\var = 3\var,
\end{split}
\end{equation}
having used, among other arguments, (\ref{main60}) and the fact that
$j \ge 1$.

Finally, consider those $n \ge m_0+n_1+m_1$ which correspond to stopping
the \o\ of $x$ during an excursion in $\tailk^c$. Say that $n = M_j+N_{j+1}+i$ for
some $j\ge 1$ and $0 \le i<m_{j+1}$. The contribution to the Birkhoff sum of the
last excursion in $\tailk^c$ is null so, in light of (\ref{main70}),
\begin{equation}
  | \ba_n f_\var(x) | = | \ba_{M_j+N_{j+1} + i} \, f_\var(x) | \le
  | \ba_{M_j+N_{j+1}}\, f_\var(x) | < 3 \var.
\end{equation}

In conclusion, $| \ba_n f_\var(x) | < 3 \var$, for \emph{all} sufficiently large
$n$, depending on $x$. Since $f - f_\var = f\, 1_{\tailk^c} \in L^1$, the proof
of Theorem \ref{thm1} is completed in the same way as the proof of
Proposition \ref{prop2}, cf.\ (\ref{main50}) \emph{et seq}.
\qed

\begin{remark}
Let us observe here that Proposition  \ref{prop2} can be thought of
as a sub-case of Theorem \ref{thm1} with $N=1$. Indeed, bearing
in mind the definition of the partition $\{L_k\}$, it is easy
to see that the condition in Proposition \ref{prop2} can be reformulated
as follows: For all $\var>0$, there exists $K\in \N \text{ such that for all }
x \in \bigcup_{k \ge K} L_k, |f(x)-f^*|\leq \var$. Here $\bigcup_{k \ge K}
L_k$ plays the role of the set $A_\var$ and this condition is exactly
that of Theorem \ref{thm1} with $N=1$.
\end{remark}

Let us now see some examples of \dsy s and \ob s to which our results
can be applied. As already mentioned, Proposition \ref{prop2} is quite
general. Consider for instance a map $T: \R \into \R$ which preserves
an \erg\ infinite, locally finite \me\ $\mu$. A nice example of such a map
is Boole's transformation $T(x) := x-1/x$, which was shown by Boole in
1857 \cite{boole} to preserve the Lebesgue measure on $\R$ and by
Adler and Weiss more than a century later \cite{aw} to be \erg.
%  to preserve the Lebesgue measure on $\R$ in 1857 (although
%the result was not stated in these terms), and shown to be ergodic by
%Adler and Weiss \cite{aw} more than a century later.
Other interesting
examples are the quasi-lifts and finite modifications thereof studied
in \cite{lmmaps}. For all the \sy s we have mentioned, every bounded
$f: \R \into \C$ such that
\begin{equation}
  f^* := \lim_{|x| \to \infty} f(x)
\end{equation}
exists verifies the hypotheses of Proposition \ref{prop2}, and therefore
$\ba f = f^*$ almost everywhere.

\fig{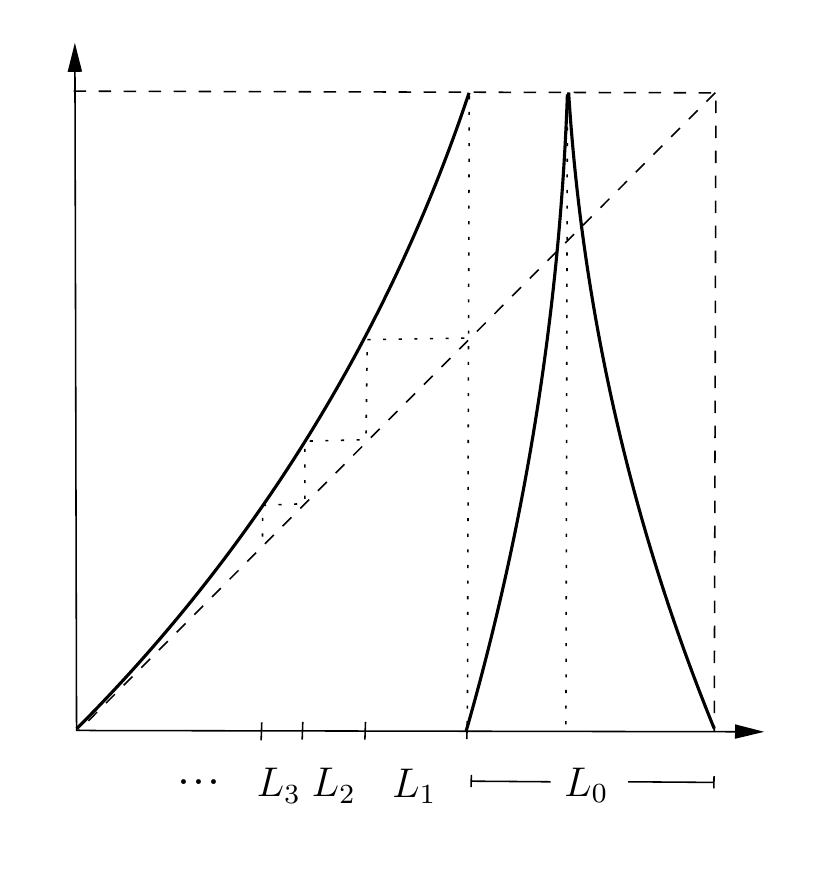}{6.7cm}{A piecewise-smooth, full-branched expanding
map of the unit interval.}

Consider now a piecewise-smooth, full-branched, expanding map
$T: [0,1] \into [0,1]$ of the type shown in Figure \ref{map01}.
If 0 is a \emph{strongly neutral} fixed point (which means that $T''$ is regular
in a neighborhood of 0), it is known \cite{t80} that under general conditions
$T$ preserves an absolutely continuous infinite \me\ $\mu$ such that
$\mu([a, 1]) < \infty$, for all $0 < a \le 1$. Also $T$ is \erg\ w.r.t.\
$\mu$ \cite{t83}. It is easy to verify that the sets $L_k$ are those marked
in Figure \ref{map01}. Therefore Proposition \ref{prop2} applies to
all bounded $f: [0,1] \into \C$ such that
\begin{equation}
  f^* := \lim_{x \to 0^+} f(x)
\end{equation}
exists.

Remaining in the case of the map $T: [0,1] \into [0,1]$ described
above, let us now introduce a class of non-trivial
global \ob s which verify the hypothesis of Theorem \ref{thm1}.
Given $N \in \Z^+$ and $c_0, c_1, \ldots, c_{N-1} \in \C$,
let $f: [0,1] \into \C$ be defined by
\begin{equation}
  f(x) = c_j \quad \Longleftrightarrow \quad x \in L_k,\: \mbox{ with } k
  \cong j \: (\mbox{mod } N).
\end{equation}
In other words, $f$ is a step \fn\ on the partition $\{ L_k \}$, which is
$N$-periodic in the index $k$. The stochastic properties of these
\ob s have been studied in \cite[Sect.\ 3.1]{bgl}. It is easy to see
that for all $x \in L_k$, with $k \ge N-1$,
\begin{equation}
  \ba_N f(x) = \frac1N \sum_{k=0}^{N-1} c_k =: f^*,
\end{equation}
Therefore Theorem \ref{thm1} applies with $N$ and $K \ge N-1$,
independently of $\eps$. Thus, $\ba f = f^*$ almost everywhere.

It is easy to extend the above idea to a class of step \fn s on
$\{ L_k \}$ which are not periodic in $k$. For example, take
a sequence $( c_k )_{k \in \N}$ of complex numbers and a number
$f^*$ such that, for every $\eps>0$, there exists $N \in \Z^+$
with the property that
\begin{equation}
  \left| \frac1N \sum_{k=j}^{j+N-1} c_k - f^* \right| \le \eps
\end{equation}
for every $j \in \N$. Examples include quasi-periodic sequences
$c_k := e^{2\pi i \alpha k}$ and many others. Then let $f$ be defined
by $f|_{L_j} \equiv c_k$. This \ob\ satisfies the hypothesis
of Theorem \ref{thm1} (again with $K \ge N-1$) by construction.

Looking for more general examples, let us consider a Kakutani tower
$\tow: Y \into Y$, as defined in (\ref{def-kt-1})-(\ref{def-kt-4}). We
also choose $L_0$ to be the base of the tower $\Sigma \times \{0\}$.
As explained in Appendix \ref{app-kakutani}, this implies that
$L_k = \rset{(x,k)} {\varphi(x) \ge k}$, that is, $L_k$ is the $k^\mathrm{th}$
level of the
tower. It is not difficult to find global \ob s $f$ which have the
approximate averaging property required by Theorem \ref{thm1} but
with different values of $f$ for different excursions outside of $L_0$.
Take for instance
\begin{equation}
  f(x,n) := e^{ 2\pi i (\omega(x) n + \gamma(x)) },
\end{equation}
defined for $(x,n) \in Y$, where $\omega$ and $\gamma$ are measurable
real-valued \fn s of $\Sigma$. Assume that there exists $\delta \in (0,1)$
such that, for all $x \in \Sigma$, $\delta \le \omega(x) \le 1-\delta$. This
implies,  for some $c = c(\delta)>0$ and for all $x \in \Sigma$, that
$| 1 - e^{-2\pi i \omega(x)} | \ge c$.

This \ob\ also satisfies the hypothesis of Theorem \ref{thm1}. In fact,
for any $\eps>0$, select $N \ge 2/c \eps$ and $K \ge N-1$. A
point $(x,n) \in \bigcup_{k \ge K} L_k$ is one for which  $n \ge K$. We
have:
\begin{equation}
  \left| \ba_N f (x,n) \right| = \frac1N \left| \sum_{k=0}^{N-1}
  e^{-2\pi i \omega(x) k} \right| \le \frac2{Nc} \le \eps.
\end{equation}
Thus $\ba f = 0$ almost everywhere.

As already alluded to in the introduction, it is hard to determine
\emph{a priori} the maximal class of global \ob s for a given \dsy.
However, for certain
\sy s, it is rather easy to agree on \fn s which ought to be considered
global \ob s; for example, \sy s defined on Euclidean spaces
(or large portions thereof) which preserve the Lebesgue \me.
In this case, the translation-invariance of the reference \me\ suggests
that at least all periodic and quasi-periodic bounded \fn s should be
global \ob s.

Let us therefore consider an interesting class of Lebesgue-\me-preserving
\dsy s in Euclidean space: piecewise-smooth, expanding maps
$T: \R^+ \into \R^+$, with full branches and an indifferent fixed point at
$+\infty$, as in Fig.~\ref{mapR+}. These maps are of the same
nature as the interval maps discussed earlier, cf.~\cite{bgl}. Indeed, if $T_o$
denotes a piecewise-smooth, full-branched,
expanding map of $[0,1]$ onto itself, and $\mu$ denotes its infinite
absolutely continuous invariant \me, then $\Phi(x) := \mu([x,1])$ defines
a bijection $(0,1) \into \R^+$ such that
$T := \Phi \circ T_o \circ \Phi^{-1}$ is a piecewise-continuous,
full-branched map $\R^+ \into \R^+$ preserving
the Lebesgue \me\ on $\R^+$. In many cases, $T$ is also
piecewise-smooth and expanding.

\fig{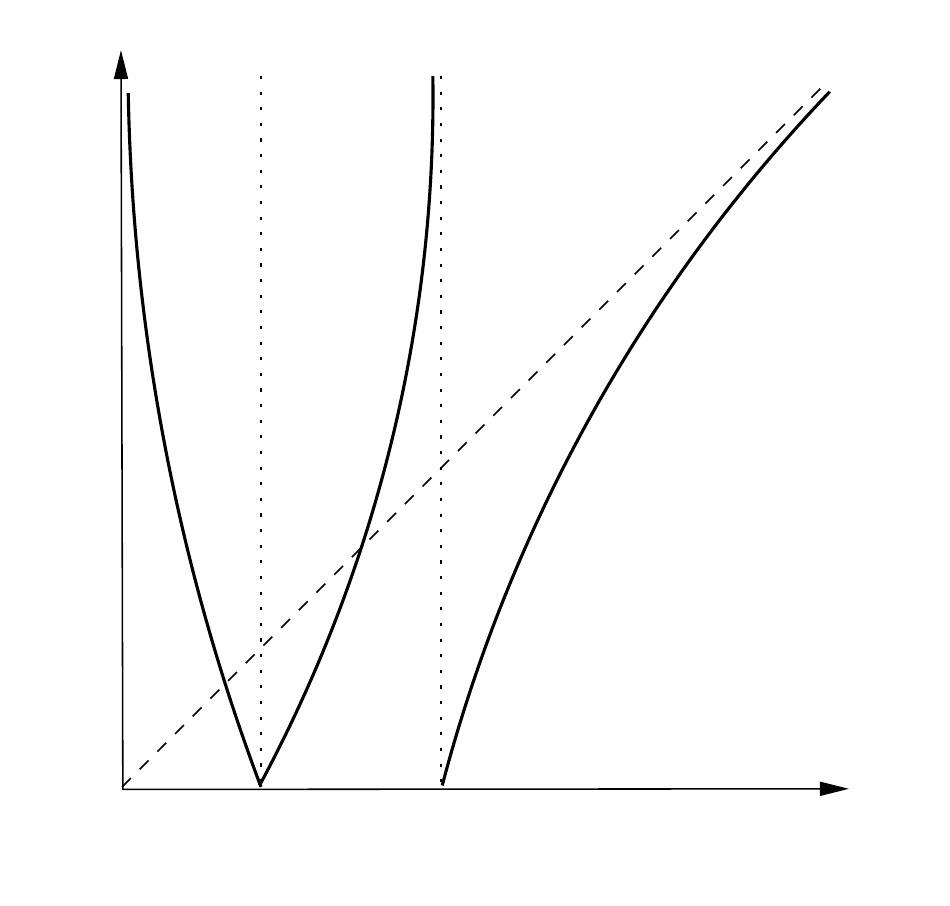} {6.7cm} {A piecewise-smooth, full-branched
expanding map of the half-line. This example corresponds to the
example of Figure \ref{map01} via the conjugation procedure
explained in the body of the paper.}

One notable example, which we will return to in Section \ref{sec-counter},
is the Farey map. This map is usually defined as a map on the unit interval,
as follows:
\begin{equation}
  F(x):=\left\{
    \begin{array}{ll}
      \ds \frac{x}{1-x}, & \mbox{for $x\in [0, 1/2]$;} \\[15pt]
      \ds \frac{1-x}{x}, & \mbox{for $x\in (1/2, 1]$.}
    \end{array}
  \right.
\end{equation}
Up to factors, $F$ has a unique Lebesgue-absolutely continuous invariant
measure $\mu$, which is given by the density $d\mu/d\leb(x) = 1/x$. Thus,
here the function $\Phi$ is given by
\begin{equation}
  \Phi(x )= \mu([x, 1]) = \int_x^1 \frac1\xi \, d\xi = -\ln x.
\end{equation}
The version of the Farey map transported to the positive real line is then
given by
\begin{equation} \label{farey-r}
  T_F(x) := -\ln(F(e^{-x})) = | \ln (e^x-1) |.
\end{equation}

Let us return to the general case of a map $T: \R^+ \into \R^+$
preserving the Lebesgue \me. The considerations made above suggest
that the first examples of global \ob s one should study are the \fn s
$f(x) := e^{2\pi i \omega x}$, with $\omega \in \R \setminus \{0\}$. Any
reasonable notion of average \cite{limix, bgl} for these \fn s would suggest that
$f^* = 0$. So the problem is to show that, for a.e.\ $x \in \R^+$, $\ba f(x) = 0$.

%In the next section, we present a class of maps for which we can show exactly this. Other examples, including that of the Farey map introduced above, we discuss in Section \ref{sec-counter}.

%Doing so would immediately yield a pointwise Birkhoff Theorem
%for functions
%\begin{equation}
%  f(x) = \int_\R e^{2\pi i \omega x} \, \sigma (d\omega),
%\end{equation}
%where $\sigma$ is a finite signed \me\ on $\R$.
%
%\nota{say what this class of \fn s is and prove that $\ba f = \sigma(0)$,
%almost everywhere.}

%\nota{To Sara: the way this came out of my mind (at this late hour),
%didn't include an introduction of the Farey map $\R^+ \into \R^+$. Not
%sure if I like it like this -- in which case Farey should be introduced together
%with its simulations, in the last section -- or the way we decided together.
%You decide. In any case, the most compact way to introduce it
%is $F(x) := | \ln (e^x-1) |$. }

\section{Periodic observables and the $\alpha$-Farey maps}
\label{sec-alpha}

In this section we present non-trivial examples of piecewise-smooth,
Lebesgue-measure-preserving, expanding maps on $\R^+$, with an indifferent fixed
point at $+\infty$, for which the Birkhoff average of $f(x) := e^{2\pi i \omega x}$,
$\omega \ne 0$, is almost everywhere zero. This will require more work than
a simple application of Theorem \ref{thm1}, but the underlying ideas are the
same.

Our maps will be conjugates, over the space $\R^+$, of the well-known
$\alpha$-Farey maps on $[0,1]$, and will be obtained by means of the construction
explained at the end of Section \ref{sec:mainthm}.

Let us recall the definition of an $\alpha$-Farey map, as introduced in
\cite[\S1.4]{kms}. Start with a decreasing sequence $(t_k)_{k \in \Z^+}$ of real
numbers such that $t_1=1$ and $\lim_{k \to \infty}t_n=0$. This sequence allows
us to define a partition
\begin{equation} \label{alpha}
  \alpha := \rset{A_k := (t_{k+1}, t_k]} {k \ge 1}
\end{equation}
of $(0,1]$. We will write $a_k := \leb(A_k) = t_k - t_{k+1}$ for the Lebesgue
measure of the $k^\mathrm{th}$ partition element. Then the map
$F_\alpha: [0,1] \into [0, 1]$ is defined by setting
\begin{equation}
  F_{\alpha}(x) := \left\{
    \begin{array}{ll}
      (1-x)/a_1, & \mbox{ for $x \in A_1$;} \\[2pt]
      {a_{k-1}} (x - t_{k+1})/a_k + t_k, & \mbox{ for $x \in A_k$, for $k \ge 2$};
      \\[2pt]
      0, & \mbox{ for $x=0$. }
    \end{array}
  \right.
\end{equation}
The map $F_\alpha$ preserves the (unique up to factors) Lebesgue-absolutely
continuous measure $\mu_\alpha$ given by the density
\begin{equation}
  h_{\alpha} := \frac{d\mu_\alpha}{d\leb} = \sum_{k=1}^\infty \frac{t_k}{a_k} \,
  1_{A_k},
\end{equation}
and the measure is infinite if and only if $\sum_k t_k = \infty$.

For later use, let us also recall the definition of the related 
$\alpha$-L\"uroth expansion (for more details we refer again to 
\cite[\S1.4]{kms}). Each partition $\alpha$ generates a series expansion of 
the numbers in the unit interval, in that we can associate to
each $x$ a sequence of positive integers $(\ell_i)_{i\geq1}$ for which
\begin{equation}
  x = t_{\ell_1} + \sum_{k=2}^\infty (-1)^{k-1} \left(\textstyle\prod\limits_{i<k}
  a_{\ell_i} \right) t_{\ell_k} = t_{\ell_1} - a_{\ell_1} t_{\ell_2} + a_{\ell_1}
  a_{\ell_2}t_{\ell_3} -\ldots
\end{equation}
To lighten the notation, we will write $x = [\ell_1, \ell_2, \ell_3, \ldots]_\alpha$.
Observe that the map $F_\alpha$ acts on this expansion in the following way:
\begin{equation}
  F_\alpha([\ell_1, \ell_2, \ell_3,\ldots]_\alpha) = \left\{
  \begin{array}{ll}
    [\ell_1-1,\ell_2, \ell_3, \ldots]_\alpha , &  \mbox{ for $\ell_1 \ge 2$}; \\[2pt]
    [\ell_2, \ell_3, \ldots]_\alpha , & \mbox{ for $\ell_1=1$} .
  \end{array}
  \right.
\end{equation}

Throughout this section, we will restrict ourselves to the particular case
$t_k := k^{-\beta}$, with $0<\beta<1/2$. The partition generated by
this sequence will be denoted by $\alpha(\beta)$. In \cite{kms} it is
referred to as an expansive partition with exponent $\beta$. Set
\begin{equation} \label{sec3-20}
  \tau_k := \sum_{j=1}^k t_j \sim \frac{ k^{1-\beta} } {1-\beta},
\end{equation}
where we write $x_k \sim y_k$ to mean that $\lim_{k \to \infty}
(x_k / y_k) =1$.

As anticipated, we want to
consider the map $T_\beta := \Phi \circ F_{\alpha(\beta)} \circ \Phi^{-1}$,
where $\Phi(x) := \mu_{\alpha(\beta)}([x, 1])$ is defined for $x \in (0,1)$.
Then $T_\beta$ is a piecewise-continuous, full-branched,
Lebesgue-preserving map on $\R^+$. In this case, $\Phi$ is a piecewise-linear
\fn\ that maps, for each $k \ge 0$, the partition element $A_{k+1}$ onto the
interval $L_k := [\tau_k, \tau_{k+1})$. Setting $L_0:=[0, \tau_1)$, a series of straightforward calculations
show that
\begin{equation}
  T_\beta(x) = \left\{
  \begin{array}{ll}
      \Phi(x), & \mbox{ for $x \in L_0$}; \\[5pt]
      \ds \frac{t_k}{t_{k+1}} (x - \tau_k) + \tau_{k-1}, & \hbox{ for $x \in L_k$,
        with $k \ge 1$}.
    \end{array}
    \right.
\end{equation}
See Fig.~\ref{aFareyR} for a picture of $T_\beta$.
Note here that, by contruction, the partition $\{L_k\}_{k \ge 0}$ matches
for $T_\beta$ the definition of the corresponding sequence of sets given
for a general map $T$ in Section \ref{sec:mainthm}.

\fig{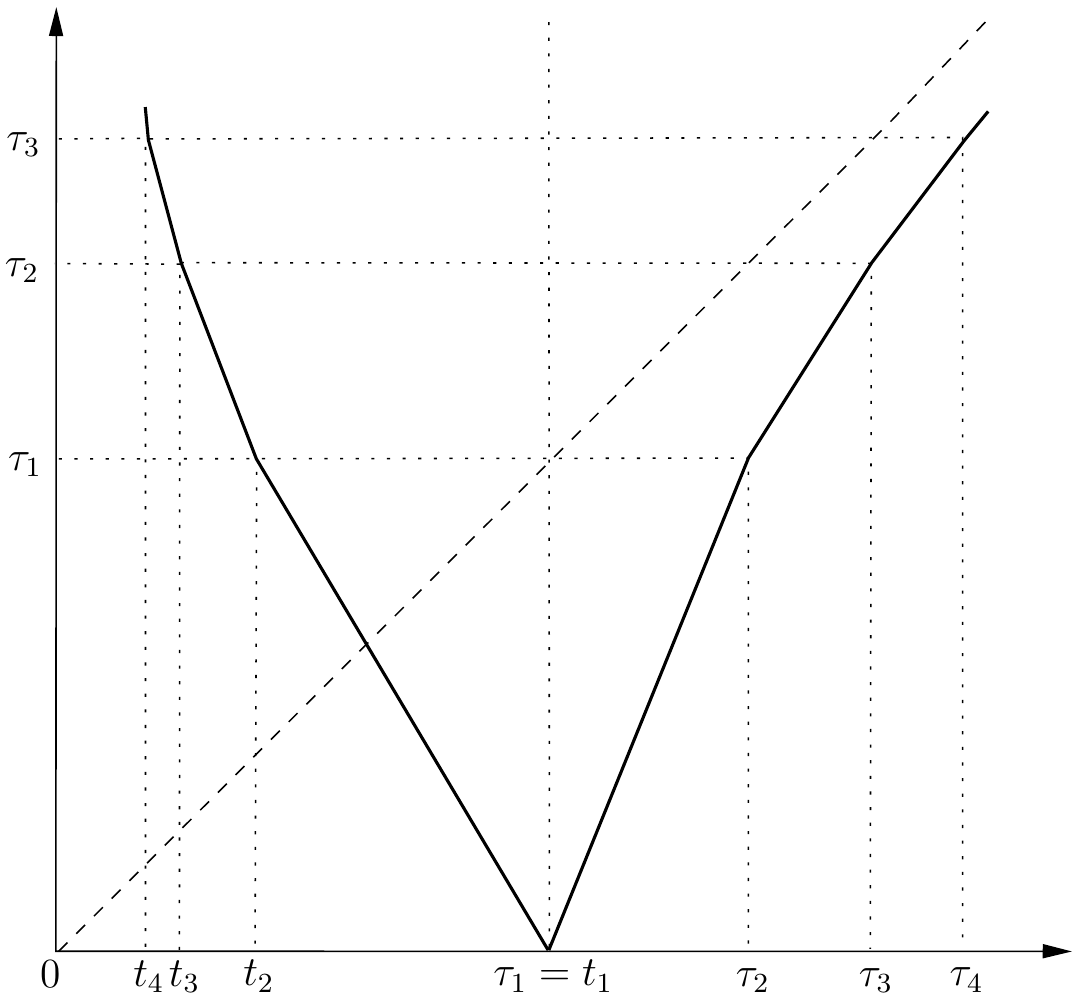} {9.6cm} {The $\alpha$-Farey map $T_\beta$ on $\R^+$.}

We are now in a position to state our main result of this section.

\begin{proposition} \label{prop3}
  Let $f(x) := e^{2\pi i \omega x}$, with $\omega \in \R \setminus \{0\}$, and
  fix $\beta \in (0,1/2)$. Then, for $\leb$-a.e.\ $x \in \R^+$,
  \begin{displaymath}
    \ba f(x) := \lim_{n\to\infty} \, \frac1n \, \sum_{i=0}^{n-1} f\circ T_\beta^i (x)
    = 0.
  \end{displaymath}
\end{proposition}

\proof First of all, let us assume without loss of generality that
$\omega>0$. Then let us define a new \ob\ $g$ which is constructed from
$f$ in the following way: For each $j \ge 0$, denote by  $k_j$ the natural
number such that $j / \omega \in L_{k_j}$. Then let
\begin{align}
  I_j &:= \bigcup_{i=k_j}^{k_{j+1}-1} \! L_i = \left[ \tau_{k_j} , \tau_{k_{j+1}}
  \right) ; \label{i-j} \\
  \omega_j &:= \frac1 {\leb(I_j)} = \frac1 { \tau_{k_{j+1}} - \tau_{k_j} }.
  \label{omega-j}
\end{align}
Finally, for all $x \in I_j$, set $g(x) := e^{2\pi i \omega_j (x - \tau_{k_j} )}$.

Let us fix one more piece of notation that we shall use throughout
the proof. For any observable $\phi$, we shall call any interval $[a,b]$
with the property that $\phi(x) = e^{2\pi i (x-a) / (b-a)}$, for $a \le x
\le b$, a {\em wavelength} for $\phi$. Therefore $g$ is a modification
of our original observable $f$ so that the wavelengths of $g$ are
unions of intervals $L_k$. Note that, since $\leb(L_k) \to 0$ when
$k \to \infty$, the modification is smaller and smaller for larger and
larger values of the argument $x$. In other words,
$\omega_k \sim \omega$.

For an arbitrary $\eps>0$, choose $K = K(\eps) \in \N$ sufficiently large that:
\begin{itemize}
\item $| f(x) - g(x) | \le \eps$ for all $x \in \tailk := \bigcup_{k \ge K} L_k$;
\item $\leb(L_K) \le \eps$, whence $\leb(L_k) \le \eps$ for all $k \ge K$;
\item $K^{\beta-1} \le \eps$.
\end{itemize}
To simplify the argument below, let us also suppose that $K=k_{j_o}$ for
some $j_o \in \N$. Moreover, $K$ will satisfy another condition which we
will state later, when it is needed.

Define the \ob\ $g_\eps := g \, 1_{\tailk}$ and consider the portion of $g_\eps$
defined on $I_j$ with $k_j \ge K$, cf.\ (\ref{i-j}). Denote by $r_j :=
k_{j+1} - k_j$ the number of intervals $L_k$ that make up the interval
$I_j$ (recall that $g_\eps$ is defined so as to have its wavelengths start and
end exactly at the endpoints of these partition elements).

Consider now a point $x \in \R^+$ whose forward \o\ intersects some
$L_k$ with $k \ge k_{j+1} - 1$. This is equivalent to asking that, at some time
$s$, $T_\beta^s(x) \in L_{k_{j+1} - 1}$. Therefore $T_\beta^{s+1}(x) \in
L_{k_{j+1} - 2}$ and so on, until  $T_\beta^{s+r_j-1} (x) \in L_{k_j}$.
In other words, the intervals $\{ L_{k_j+i} \}_{i=0}^{r_j-1}$ that partition $I_j$
each contain exactly one \o\ point of $x$, from time $s$ to time
$s+r_j-1$. We want to compare these intervals to the intervals
$\{ B_{k_j+i} \}_{i=0}^{r_j-1}$, which are defined to be a partition of $I_j$
into intervals of the same size, labeled from left to right.

We claim that, for $0 \le i < r_j$,
\begin{equation}\label{sec3-10}
  L_{k_j+i} \cap B_{k_j+i} \ne \emptyset.
\end{equation}
Indeed, observe first that the common size of the intervals
$B_{k_j+i}$ is the average of the sizes of the intervals
$L_{k_j+i}$, therefore the ``relative discrepancy'' between the
sizes of corresponding sets can be estimated as follows:
\begin{equation}\label{sec3-12}
  \frac{ |m(L_{k_j+i}) - m(B_{k_j+i})| } {m(B_{k_j+i})} = \left| \frac{
  \leb( L_{k_j+i} ) } { \leb( B_{k_j+i} ) } - 1 \right| \le \frac{ \leb( L_{k_j} ) }
  { \leb( L_{k_{j+1}} ) } -1,
\end{equation}
because $\leb( L_k )$ is a decreasing \fn\ of $k$. Clearly, the relative
discrepancy between the sizes of $\bigcup_{i=0}^q L_{k_j+i}$ and
$\bigcup_{i=0}^q B_{k_j+i}$, for $0 \le q < r_j$, does not exceed
the sum of the individual discrepancies (\ref{sec3-12}). Condition
(\ref{sec3-10}) will be satisfied if the former is always less than or
equal to 1. A sufficient condition for this is
\begin{equation}\label{sec3-13}
  r_j \left(\frac{ \leb( L_{k_j} ) } { \leb( L_{k_{j+1}} ) }-1\right) \le 1.
\end{equation}

Towards the proof of (\ref{sec3-13}), we make
several observations. First, by construction, $\leb(L_k) = t_{k+1} \sim
k^{-\beta}$. On the other hand, the definition of $k_j$ and (\ref{sec3-20})
imply that, as $j \to \infty$,
\begin{equation} \label{sec3-14}
  \frac{j}{\omega} \sim \tau_{k_j} \sim \frac{ k_j^{1-\beta} } {1-\beta},
\end{equation}
whence
\begin{equation} \label{sec3-15}
  k_j \sim \left( \frac{1-\beta}{\omega} \, j\right)^{1/(1-\beta)}
\end{equation}
and
\begin{equation} \label{sec3-80}
  \frac{ \leb( L_{k_j} ) } { \leb( L_{k_{j+1}} ) } \sim \left( 1 + \frac1j
  \right)^{\beta/(1-\beta)} \!\! \sim  1+ \frac{\beta}{1-\beta} \, \frac1j.
\end{equation}
Also, by (\ref{sec3-15}),
\begin{equation}\label{sec3-90}
  r_j := k_{j+1}-k_j \sim c_1 \, j^{\beta / (1-\beta)},
\end{equation}
for some $c_1 = c_1(\beta, \omega) > 0$.

Putting these observations together, we obtain that, for some positive
constant $c_2$,
\begin{equation} \label{sec3-95}
  r_j \left(\frac{ \leb( L_{k_j} ) } { \leb( L_{k_{j+1}} ) } -1 \right) \sim c_2 \,
  j^{(2\beta-1)/(1-\beta)} .
\end{equation}
Since $0 < \beta < 1/2$, the above term vanishes as
$j \to \infty$. If we choose $j_o$ sufficiently large, that is, if we choose
$K = k_{j_o}$ sufficiently large (which is the condition we anticipated
earlier we would state precisely), we can guarantee that (\ref{sec3-13})
holds for all $j \ge j_o$. This proves the claim (\ref{sec3-10}).

Now, for $k \in \N$, denote by $b_k$ the midpoint of the interval $B_k$.
Recalling that $\leb(L_k) < \var$ for all $k\ge K$ and that
$T_\beta^{s+i}(x) \in L_{k_{j+1} - 1 - i}$, for all $0 \le i < r_j$, it follows
from (\ref{sec3-10}) that, for the same values of $i$,
\begin{equation} \label{sec3-97}
  \left| T_\beta^{s+i}(x) - b_{k_{j+1} - 1 - i} \right| < 2\var.
\end{equation}
Then, since $g_\var$ is Lipschitz continuous on its wavelength $I_j$, with
constant $2\pi \omega_j$, and we can find an upper bound $c_3>0$ such
that $2\pi \omega_j \le c_3$ for all $j \ge j_o$, we have
\begin{equation} \label{sec3-98}
  \left| g_\var( T_\beta^{s+i}(x) ) - g_\var( b_{k_{j+1} - 1 - i} ) \right| <
  2\, c_3\, \var.
\end{equation}
On the other hand,
\begin{equation} \label{sec3-99}
  \sum_{i=0}^{r_j-1} g_\var ( b_{k_{j+1} - 1 - i} ) = 0
\end{equation}
because, for $k_j \le k < k_{j+1}$, the $b_k$ are the midpoints of the
uniform partition of $I_j$. Therefore
\begin{equation} \label{sec3-100}
  \left| \frac1{r_j} \sum_{i=0}^{r_j-1} g_\var (T_\beta^{s+i}(x)) \right| <
  2\, c_3\, \var.
\end{equation}

In summary, if we have a section of orbit of a point $x$ under $T_\beta$
that travels through an entire wavelength of the function $g_\var$, then the
partial Birkhoff average through this excursion can be at most $2 c_3 \var$.

Now, for $x \in \R^+$, let $(n_i)_{i \in \Z^+}$ denote the sequence of hitting
times to $L_0$. More precisely, $x \in L_{n_1}$, $T_\beta^{n_1}(x) \in L_0$,
$T_\beta^{n_1+1}(x) \in L_{n_2}$, and so on. (Note that these digits are
closely related to the $\alpha(\beta)$-L\"uroth digits of $\Phi^{-1}(x) \in [0, 1]$,
which are $[n_1+1, n_2+1, n_3 +1, \ldots]_{\alpha(\beta)}$.)

We want first to consider $\ba_n g_\var$, with $n =
N_q := \sum_{i=1}^q n_i$. In other words, we want to estimate
the Birkhoff average of $g_\var$ along an entire number of excursions
back to $L_0$. Let us suppose first of all that $n_i \ge K$ for all $i \ge 1$,
since  otherwise we would be adding only zeros for certain portions of the
orbit. Consider the portion of the orbit that lies in the sets $L_{n_i}, L_{n_i-1},
\ldots, L_K$. This can be split into $p_i +1$ blocks, where the
initial block runs through a portion (in general) of a wavelength of
$g_\var$ and the other $p_i$ blocks run through complete
wavelengths. Using the index $0 \le u \le p_i$, to denote these blocks,
where $u=0$ refers to the first block, which corresponds to the
partial wavelength, let $\rho^{(i)}_u$ be the number of intervals $L_k$ in
the $u^\mathrm{th}$ block. In other words, if the $u^\mathrm{th}$ block
corresponds to the wavelength $I_j$, then $\rho^{(i)}_u = r_j$. Notice
that, by (\ref{sec3-14}) and (\ref{sec3-90}), $r_j$ is asymptotic to
$k_j^\beta \le n_i^\beta$. This shows in particular that $\rho^{(i)}_0 \le
c_4 n_i^\beta$, for some constant $c_4 = c_4(\beta, \omega)$.

So, in the special case $n = N_q = \sum_{i=1}^q ( \sum_{u=0}^{p_i}
\rho^{(i)}_u + K)$, we obtain
\begin{equation} \label{sec3-110}
\begin{split}
  \left| \ba_n g_\var(x) \right| &< \frac1{N_q} \sum_{i=1}^q \left( c_4\,
  n_i^\beta + 2\, c_3 \, \var \sum_{u=1}^{p_i} \rho^{(i)}_u \right) \\
  &< c_4 \, K^{\beta-1} \frac1 {N_q} \sum_{i=1}^q n_i  +
  2\, c_3 \, \var \\[4pt]
  &\le \left( c_4 + 2\, c_3 \right) \var =: c_5 \, \var.
\end{split}
\end{equation}
Notice that in the second inequality we have used the fact that
$n_i^\beta = n_i^{\beta-1}n_i\leq K^{\beta-1}n_i$ and the estimate
(\ref{sec3-100}). In the third inequality we have used the assumption
$K^{\beta-1} \le \var$.

Let us now consider the case $N_q < n < N_{q+1}$,
i.e., we consider the Birkhoff average of $g_\var$ up to a point
which is in the middle of an excursion back to $L_0$. In the portion of \o\
between time $N_q$ and time $n$, there could be up to two
blocks (an initial and a final block) that are neither contained in
$\bigcup_{k=0}^{K-1 }L_k$ or form a full wavelength. The contribution
to the Birkhoff sum from these blocks, which is of order at most
$n_{q+1}^\beta$, cannot cannot be compensated for as in
(\ref{sec3-110}), because the denominator $n > N_q$ might be
much smaller than $N_q + n_{q+1} = N_{q+1}$. This is a phenomenon
that occurs because the numbers $(n_i)$ are distributed (in a sense
better specified in the proof of Lemma \ref{lem:errors} below) as the
outcomes of a non-integrable random variable (more precisely, a
random variable in the domain of attraction of a $\beta$-stable
distribution). Hence, for some $q$, the number $n_{q+1}$ might be very
large compared to $N_q$.

We appeal instead to the following lemma, which we prove at the end
of this section.

\begin{lemma} \label{lem:errors}
  Let $x = [\ell_1, \ell_2, \ldots]_{\alpha(\beta)} \in [0, 1]$ denote the
  $\alpha(\beta)$-L\"uroth expansion of the point $x$, where $0<\beta<1$.
  Then, for Lebesgue-almost every $x\in [0, 1]$,
  \begin{displaymath}
    \lim_{n\to\infty} \, \frac{\ell_n^\beta}{\ell_1+\cdots+\ell_{n-1}}=0.
  \end{displaymath}
\end{lemma}

Since the conjugation $\Phi: (0,1) \into \R^+$ is non-singular and
since, for a general $n$, $| \ba_n g_\var(x) |$ does not exceed
a constant times
\begin{equation}
  \frac{\sum_{i=1}^{q+1} n_i^\beta} {\sum_{i=1}^q n_i} +
  2\, c_3\, \var \, \frac{\sum_{i=1}^{q+1} \sum_{u=1}^{p_i} \rho^{(i)}_u}
  {\sum_{i=1}^{q+1} \sum_{u=0}^{p_i} \rho^{(i)}_u},
\end{equation}
we conclude that, in view of Lemma \ref{lem:errors} and estimate
(\ref{sec3-110}),
\begin{equation}
  \limsup_{n \to \infty} \left| \ba_n g_\var(x) \right| \le c_5\, \var
\end{equation}
for a.e.\ $x \in \R^+$. Defining $f_\eps := f \, 1_{\tailk}$ and recalling
that, by the assumption on $K$, $\| f_\eps - g_\eps \|_\infty \le \eps$, we
deduce that, for all $\eps>0$, there exists a Lebesgue-full-\me\ set
$B_\eps \subseteq \R^+$ such that, for every $x \in B_\eps$,
$\limsup_{n \to \infty} | \ba_n f_\eps(x) | \le (c_5+1) \eps$. Since
$f - f_\eps \in L^1(\R^+, \leb)$, the conclusion of Proposition \ref{prop3}
is achieved in the same way as that of Proposition \ref{prop2}.
\qed

\proofof{Lemma \ref{lem:errors}} We first claim that to obtain the
statement of the lemma it is enough to show the following: For all $\var>0$,
\begin{equation} \label{err-10}
  \sum_{n=2}^\infty \leb\! \left( \ell_n^\beta \ge \var \sum_{i=1}^{n-1} \ell_i
  \right)<\infty.
\end{equation}
Indeed, if this relation holds, by the Borel-Cantelli lemma we infer that
\begin{equation}
  \leb \! \left( \rset{ x = [\ell_1, \ell_2, \ldots]_{\alpha(\beta)}\in [0,1] }
  {\ell_n^\beta\geq \var \sum_{i=1}^{n-1}\ell_i \mbox{ for infinitely many }
  n \in \Z^+} \right) = 0.
\end{equation}
In other words, there exists $B_\eps \subseteq [0,1]$, with $\leb( B_\eps )
=1$, such that
\begin{equation}
  \frac{\ell_n^\beta} {\sum_{i=1}^{n-1} \ell_i} \le \var,
\end{equation}
for all $n$ larger than some $N = N(x, \eps)$. Fix a vanishing sequence
$( \eps_i)_{i \in \N}$ and define $B := \bigcap_i B_{\eps_i}$. Clearly
$\leb(B)=1$ and, for all $x \in [\ell_1, \ell_2, \ldots]_{\alpha(\beta)} \in B$,
the limit in the statement of the lemma holds true.

It thus remains to prove (\ref{err-10}). We start by observing that for every
$\alpha$-L\"uroth map the digits $(\ell_i)$ are independent identically-distributed
random variables w.r.t.\ $\leb$. We then make the following
sequence of observations:
\begin{equation} \label{err-20}
\begin{split}
  \sum_{n=2}^\infty \leb\! \left(\ell_n^\beta \ge \var \sum_{i=1}^{n-1} \ell_i
  \right) &= \sum_{n=2}^\infty \sum_{k={n-1}}^\infty \leb\! \left( \left(
  \ell_n^\beta \ge \var k \right) \cap \left( \sum_{i=1}^{n-1} \ell_i = k
  \right) \right) \\
  &= \sum_{n=2}^\infty \sum_{k={n-1}}^\infty \leb\! \left( \ell_n \ge
  (\var k)^{1/\beta} \right) \leb\! \left( \sum_{i=1}^{n-1} \ell_i = k \right) \\
  &= \sum_{k=1}^\infty \sum_{n=2}^{k+1} \leb\! \left( \ell_n \ge
  (\var k)^{1/\beta} \right) \leb\! \left( \sum_{i=1}^{n-1}\ell_i = k \right).
\end{split}
\end{equation}
Now, for $k \in \Z^+$, let us define
\begin{equation}
  \mathcal{C}_k^{(\alpha(\beta))} := \rset{ x=[\ell_1, \ell_2, \ldots]_{\alpha} }
  { \exists n \mbox{ with } \sum_{i=1}^n \ell_i=k }.
\end{equation}
In the language of \cite{kms} this is a \emph{sum-level set} for the
partition $\alpha(\beta)$, or an $\alpha(\beta)$-sum-level set. Observe
that, when such an $n$ exists, clearly $n\le k$. So, in light of the fact that
$(\ell_n)$ are i.i.d., we can then rewrite (\ref{err-20}) as
\begin{equation}
  \sum_{n=2}^\infty \leb\! \left (\ell_n^\beta \ge \var \sum_{i=1}^{n-1}\ell_i \right)
  = \sum_{k=1}^\infty \leb\! \left( \ell_1\ge (\var k)^{1/\beta} \right) \leb\! \left(
  \mathcal{C}^{(\alpha(\beta))}_k \right).
\end{equation}
We have already mentioned that $\alpha(\beta)$ is an expansive partition
with exponent $\beta$, therefore, by Theorem 1(2)(ii) of \cite{kms},
\begin{equation}
  \leb\! \left( \mathcal{C}^{(\alpha(\beta))}_k \right) \sim \frac1
  {\Gamma(2-\beta) \Gamma(\beta)}  \left( \sum_{j=1}^k t_j \right)^{-1}
  \!\! \sim \frac{c} {k^{1-\beta}},
\end{equation}
where $\Gamma$ denotes Euler's Gamma function and $c$ is a positive
constant (recall that $t_k := k^{-\beta}$). Observing also that
\begin{equation}
  \leb\! \left(\ell_1 \ge (\var k)^{1/\beta} \right) = t_{\lceil(\var k)^{1/\beta}\rceil}
  = \left\lceil(\var k)^{1/\beta} \right\rceil^{-\beta} \sim (\var k)^{-1},
\end{equation}
we finally obtain that
\begin{equation}
  \sum_{n=2}^\infty \leb\! \left(\ell_n^\beta \ge \var \sum_{i=1}^{n-1} \ell_i
  \right) \sim \frac{c}{\var} \, \sum_{k=1}^\infty \frac{1}{k^{2-\beta}} <\infty.
\end{equation}
This gives (\ref{err-10}) and concludes the proof of Lemma \ref{lem:errors}.
\qed

\begin{remark} \label{rk-prop3-1}
  Proposition \ref{prop3} can be extended to include the
  case $\beta = 1/2$. In that case, in fact, the r.h.s.\ of (\ref{sec3-95})
  does not vanish as $j \to \infty$, but it is still bounded above.
  In other words, the inequality (\ref{sec3-13}) holds with a
  bound possibly larger than 1 on the r.h.s. The pivotal relation
  (\ref{sec3-10}) may not hold anymore, but one can still claim that
  the distance between the intervals $L_{k_j+i}$ and $B_{k_j+i}$ is
  a bounded multiple of $\leb( B_{k_j+i} )$, which tends to zero as
  $j \to \infty$. So it suffices to select a large enough $K = k_{j_o}$
  to guarantee that the l.h.s.\ of (\ref{sec3-97}) does not exceed
  $c_6 \var$, for some $c_6>0$. The rest of the proof holds, with
  possibly different constants.
\end{remark}

\begin{remark}
  Another way to generalize Proposition \ref{prop3} is by proving its
  statement for any periodic continuous global \ob\ $f$.
  In such case, one calls ``wavelength'' any interval of the type
  $[jT, (j+1)T]$, where $T$ is the period of $f$. Without loss of
  generality, as we have seen many times so far, one can assume
  that $f^* := \int_{jT}^{(j+1)T} \! f \, d\leb = 0$. As in the proof of the
  proposition, one constructs the \ob\ $g$ by means of a piecewise
  affine transformation that adapts the wavelengths of $f$ to the
  intervals $I_j$ defined by (\ref{i-j}). Once again $|f(x) - g(x)| \le \eps$,
  for all large enough $x$. The proof flows as before, except that:
  \begin{enumerate}
  \item The restrictions $g_\var |_{I_j} = g |_{I_j}$ are not Lipschitz
    continuous but only continuous (thus uniformly continuous), so
    (\ref{sec3-98}) might not hold. On the other hand, since the graphs
    of $g |_{I_j}$ get closer and closer, upon suitable translation, as
    $j \to \infty$, it is easy to see that one can find a \fn\
    $\delta \mapsto \Omega(\delta)$, with $\lim_{\delta \to 0^+}
    \Omega(\delta) = 0$, which is an upper bound for the moduli of
    continuity of all $g |_{I_j}$, for $j \ge j_o$. Therefore, (\ref{sec3-98})
    can be replaced by
    \begin{equation} \label{sec3-98-prime}
      \left| g_\var( T_\beta^{s+i}(x) ) - g_\var( b_{k_{j+1} - 1 - i} ) \right|
      < \Omega(2\eps).
    \end{equation}

  \item Equation (\ref{sec3-99}) may not be true: the Riemann sum of
    $g$ over the midpoints of the uniform partition of $I_j$ is not
    \emph{exactly} zero in general, but it is nonetheless close to zero if
    the partition is dense enough, which happens for $j$ large enough.
    In other words, for every $j \ge j_o$ (with a possible redefinition of
    $j_o$),
    \begin{equation} \label{sec3-99-prime}
      \left| \frac1{r_j} \sum_{i=0}^{r_j-1} g_\var ( b_{k_{j+1} - 1 - i} ) \right|
      \le \eps.
    \end{equation}
  \end{enumerate}
  So, replacing (\ref{sec3-98})-(\ref{sec3-99}) with
  (\ref{sec3-98-prime})-(\ref{sec3-99-prime}), one rewrites
  (\ref{sec3-100}) as
  \begin{equation} \label{sec3-100-prime}
    \left| \frac1{r_j} \sum_{i=0}^{r_j-1} g_\var (T_\beta^{s+i}(x)) \right| <
    \Omega(2\eps) + \eps.
  \end{equation}
  What matters is that the above r.h.s.\ vanishes for $\eps \to 0^+$. The
  rest of the proof is the same, except that one uses $\Omega(2\eps) +
  \eps$ instead of $2c_3 \eps$.
\end{remark}

\section{Counterexamples and discussion}
\label{sec-counter}

The results of Section \ref{sec:mainthm} are obviously not as strong
as the original Birkhoff Theorem, which covers a wide functional space
of \ob s (that is, $L^1$). In  particular, one may point out that the main
hypothesis of Theorem \ref{thm1} --- namely, in overwhelmingly large
portions of the space, the partial Birkhoff averages of $f$ over long
(though fixed!) times are well approximated by a constant --- somehow
contains the assertion of the theorem. On the other hand, the known
examples of $L^\infty$ \ob s whose Birkhoff average does not
converge almost everywhere to a constant are precisely \fn s that are
constant on regions of the space where the moving point takes
longer and longer excursions.

The most famous such example concerns Boole's transformation
$T(x) = x - 1/x$ on $\R$, cf.\ Section \ref{sec:mainthm}. It is known
that the \emph{frequency} of visits to the positive half-line follows
a non-trivial law \cite{t02}. More precisely, if $\nu$ is a
Lebesgue-absolutely continuous probability \me\ on $\R$ (remember
that $T$ preserves the Lebesgue \me) then, for $0 \le t \le 1$,
\begin{equation} \label{arcsin-law}
  \lim_{n \to \infty} \nu( \{ \ba_n (1_{\R^+}) \le t \} ) =
  \frac 2 \pi \arcsin \sqrt{t}.
\end{equation}
Thus, $\ba_n$ cannot converge to a constant almost everywhere.
As a matter of fact, it can be proved that, for $\leb$-a.e.\ $x \in \R$,
\begin{equation}
\begin{split}
  \liminf_{n \to \infty} \ba_n (1_{\R^+}) &= 0, \\
  \limsup_{n \to \infty} \ba_n (1_{\R^+}) &= 1.
\end{split}
\end{equation}
The limit in (\ref{arcsin-law}) is usually referred to as the \emph{arcsine
law} for the occupation times of half-lines for Boole's transformation.
In fact, since $1_{[a, +\infty)} - 1_{\R^+} \in L^1(\R, \leb)$ for all $a \in
\R$, it follows from Proposition \ref{prop1} that the indicator \fn\ of
$\R^+$ in (\ref{arcsin-law}) can be replaced by that of any other
right half-line. The law for $\ba_n (1_{(-\infty, a]})$ follows
straightforwardly.

This example has been generalized in a number of ways, cf.\
\cite{tz,sy} and references therein. A strong recent result is that of
Sera and Yano \cite[Thm.~2.7]{sy} about the joint distribution of
the frequencies of visits to many infinite-measure sets. We describe it in
loose terms: For an infinite-\me-preserving, conservative, \erg\
\dsy\ $(X, \mu, T)$, suppose that $X$ is the disjoint union of
$X_0, R_1, R_2, \ldots, R_d$, with $0 < \mu(X_0) < \infty$ and
$\mu(R_i) = \infty$ for all $i$. Suppose also that, for $i \ne j$, an \o\
point cannot pass from $R_i$ to $R_j$ without visiting $X_0$. One
says that the \emph{rays} $\{ R_i \}$ are \emph{dynamically separated}
by the \emph{junction} $X_0$. Under a couple of important technical
assumptions (one concerning the so-called asymptotic entrance densities
from the rays into the junction and the other concerning the regular
variation of certain normalizing rates \cite[Ass's.~2.3 and 2.5]{sy})
the random vector $( \ba_n ( 1_{R_1} ) , \ldots, \ba_n
( 1_{R_d} ) )$ converges in distribution, w.r.t.\ any probability
$\nu \ll \mu$, to the vector
\begin{equation}
  \frac{ (\xi_1, \ldots, \xi_d) } {\sum_{i=1}^d \xi_i},
\end{equation}
where $\xi_1, \ldots, \xi_d$ are positive i.i.d.\ random variables with
one-sided $\beta$-stable distributions (save for degenerate cases).
Therefore, typically, the Birkhoff average
\begin{equation} \label{ob-rays}
  f = \sum_{i-1}^d \gamma_i \, 1_{R_i},
\end{equation}
converges in distribution of a non-constant variable.

The above discussion shows that some ``dynamical averaging''
hypothesis is needed for a bounded \ob\ to fulfill the assertion of the
Birkhoff Theorem. So the question is, how slowly is a \fn\ allowed
to vary over the \o s of $T$ (say, how close to a constant on each ray
must it be) in order for it to still have a constant overall Birkhoff average?

The \dsy\ and \ob s of Section \ref{sec-alpha} are a good case
study. Proposition \ref{prop3} and Remark \ref{rk-prop3-1} guarantee
that, for $0 < \beta \le 1/2$, the Birkhoff average of the `wave'
$f(x) = e^{2\pi i \omega x}$,
for the $\alpha(\beta)$-Farey map $T_\beta$, vanishes almost
everywhere. Recall that the lengths of the cylinders $L_k$ of $T_\beta$
decrease like $k^{-\beta}$. This means that, the smaller the $\beta$,
the closer the partition $\{ L_k \}$ is to the uniform partition, when
restricted to a wavelength of $f$, implying that the \o\ segments that
traverse a wavelength of $f$ contribute with an almost null partial
average for $f$. By contrast, for $\beta$ close to 1, there will be
many more cylinders in the right half of the wavelength than in the
left half, making the variation of $f$ along an \o\ segment in the
right half much slower than the corresponding variation on the
left half.

So the arguments in the proof of Proposition \ref{prop3} do not
work for $\beta > 1/2$. We do not know whether $\ba_n f$ vanishes
for $\beta > 1/2$ as well, but numerical simulations do show a
different behavior than the case $\beta \le 1/2$; see
Figs.~\ref{beta_35}--\ref{beta_65}.

\figsimul{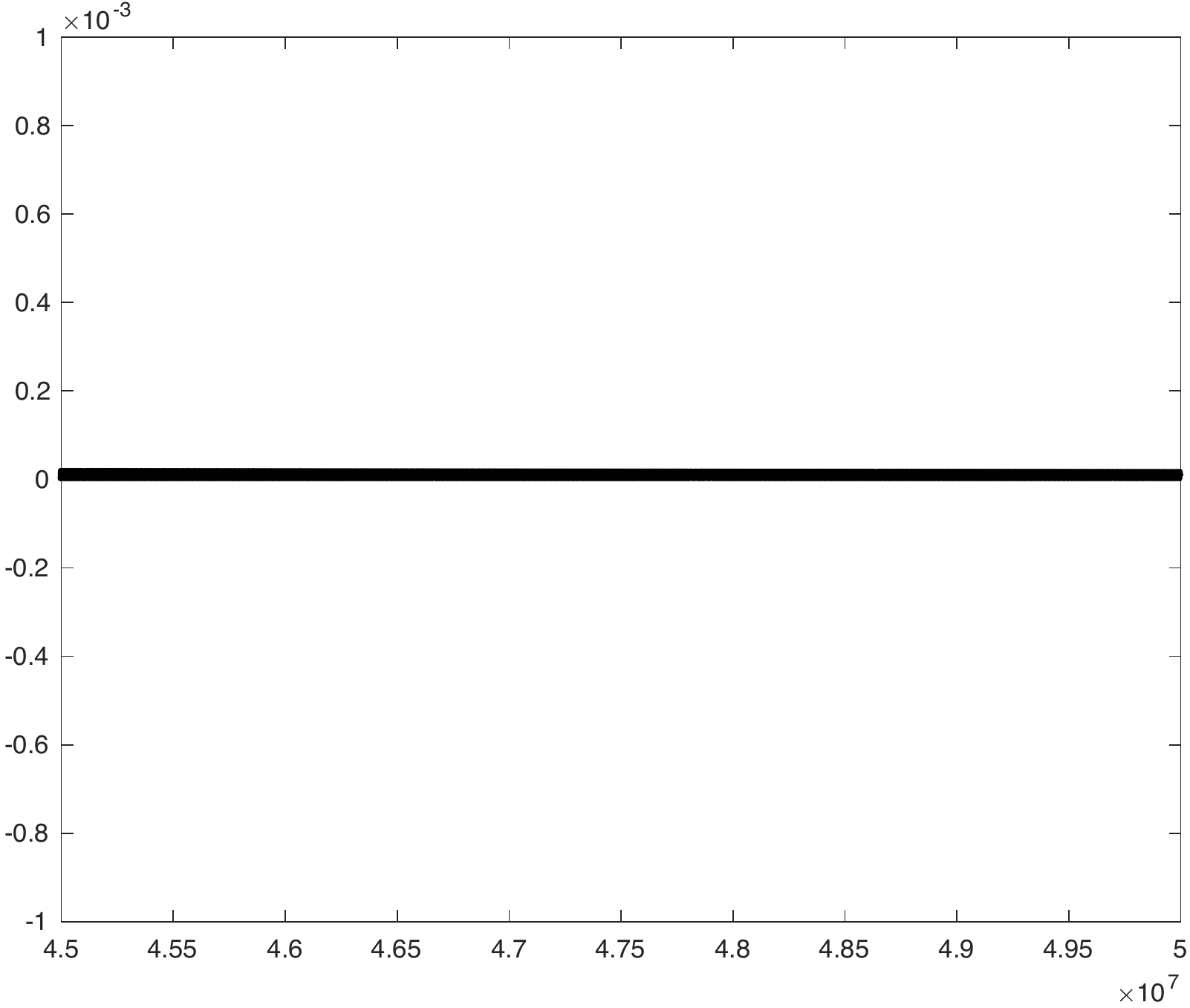}{For the map $T_\beta$ with $\beta=0.35$, the
figure shows a plot of $\ba_n g(x_0)$, with $g(x) = \cos( 2\pi \omega x)$,
$\omega=0.2$ and $x_0 = 0.65$. Here $4.5 \cdot 10^7 \le n \le 5 \cdot
10^7$ and the vertical scale is in units $10^{-3}$.}

\figsimul{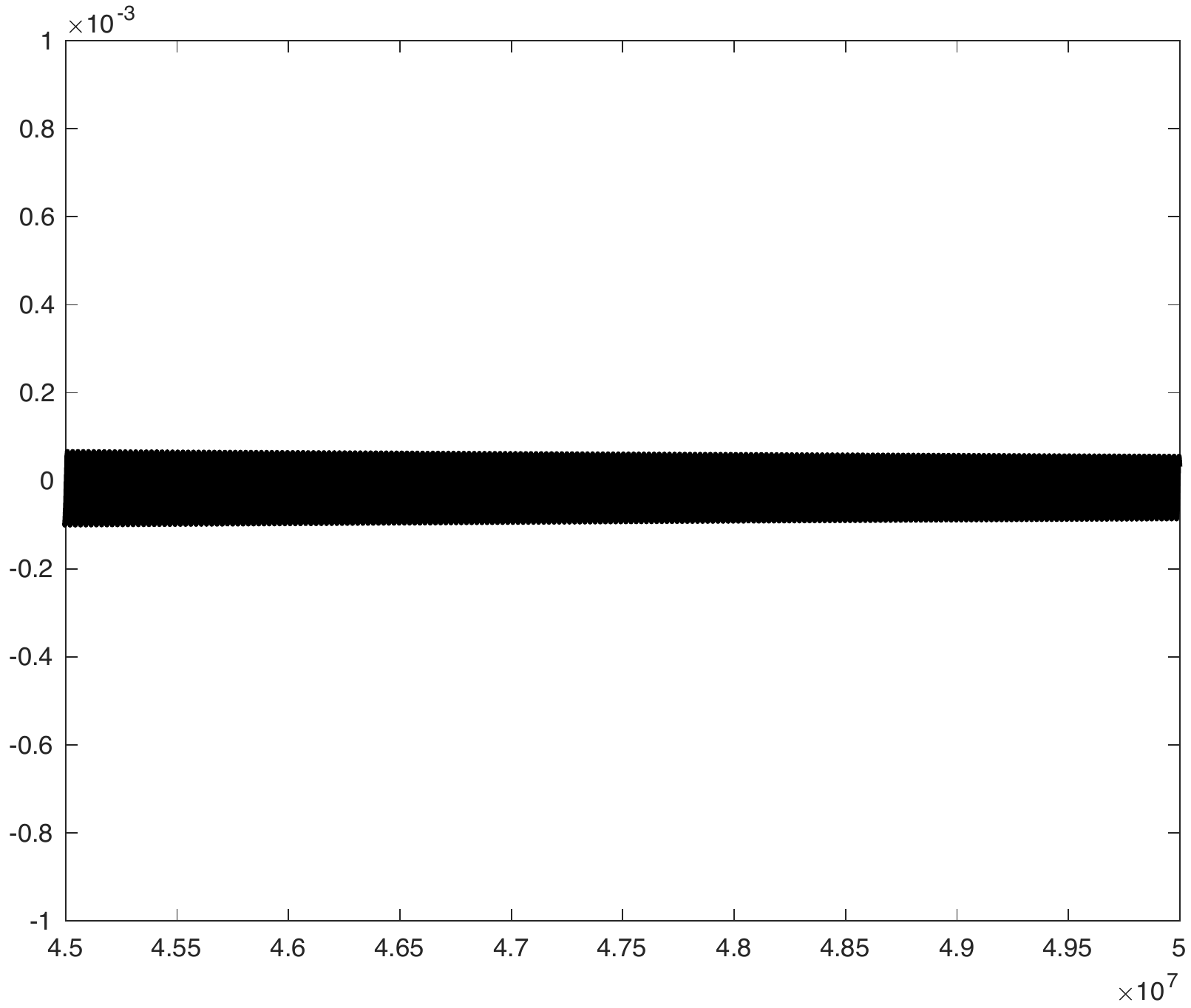}{Same plot as in Fig.~\ref{beta_35}, for the case
$\beta=0.48$.}

\figsimul{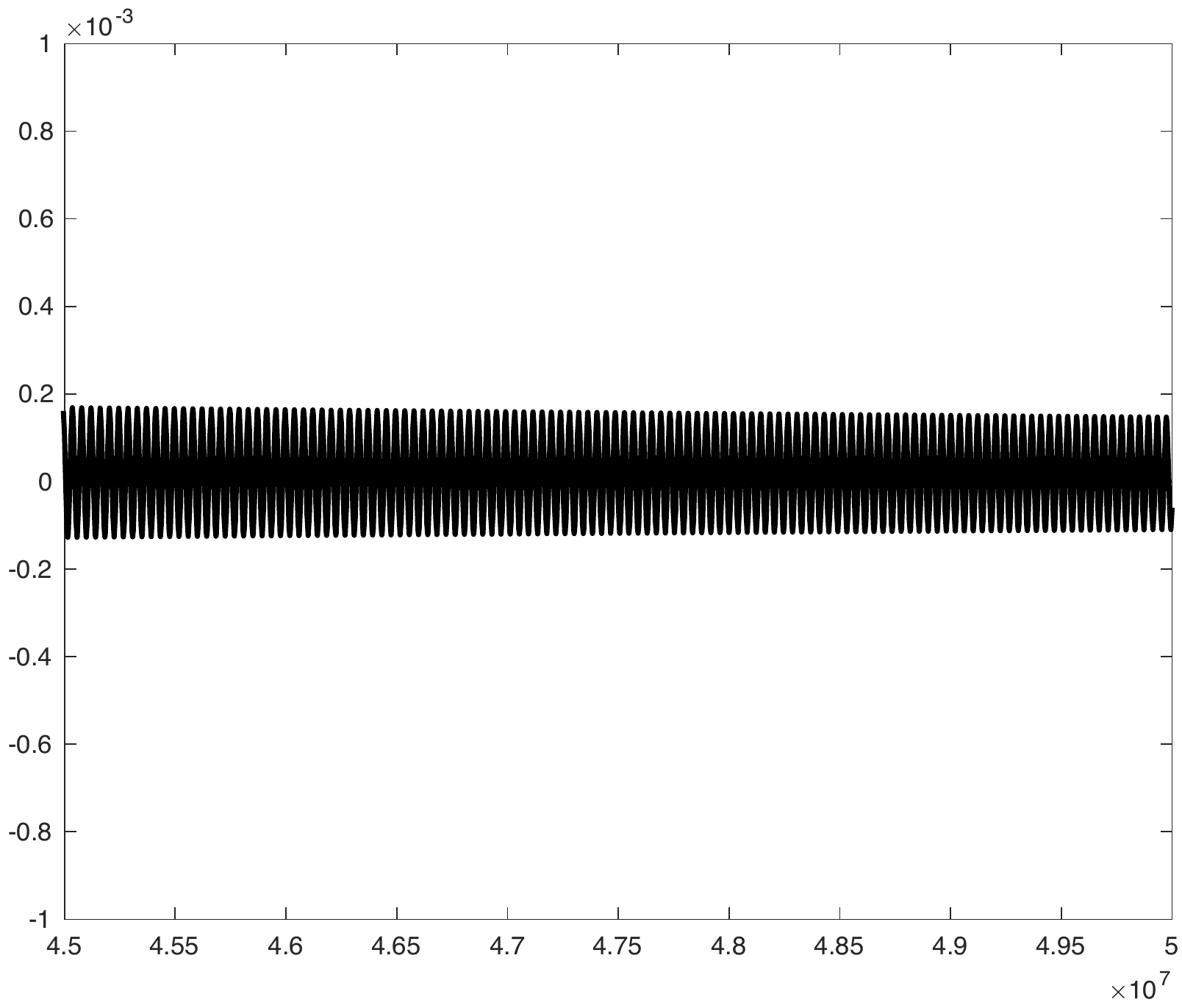}{Same plot as in Fig.~\ref{beta_35}, for the case
$\beta=0.50$.}

\figsimul{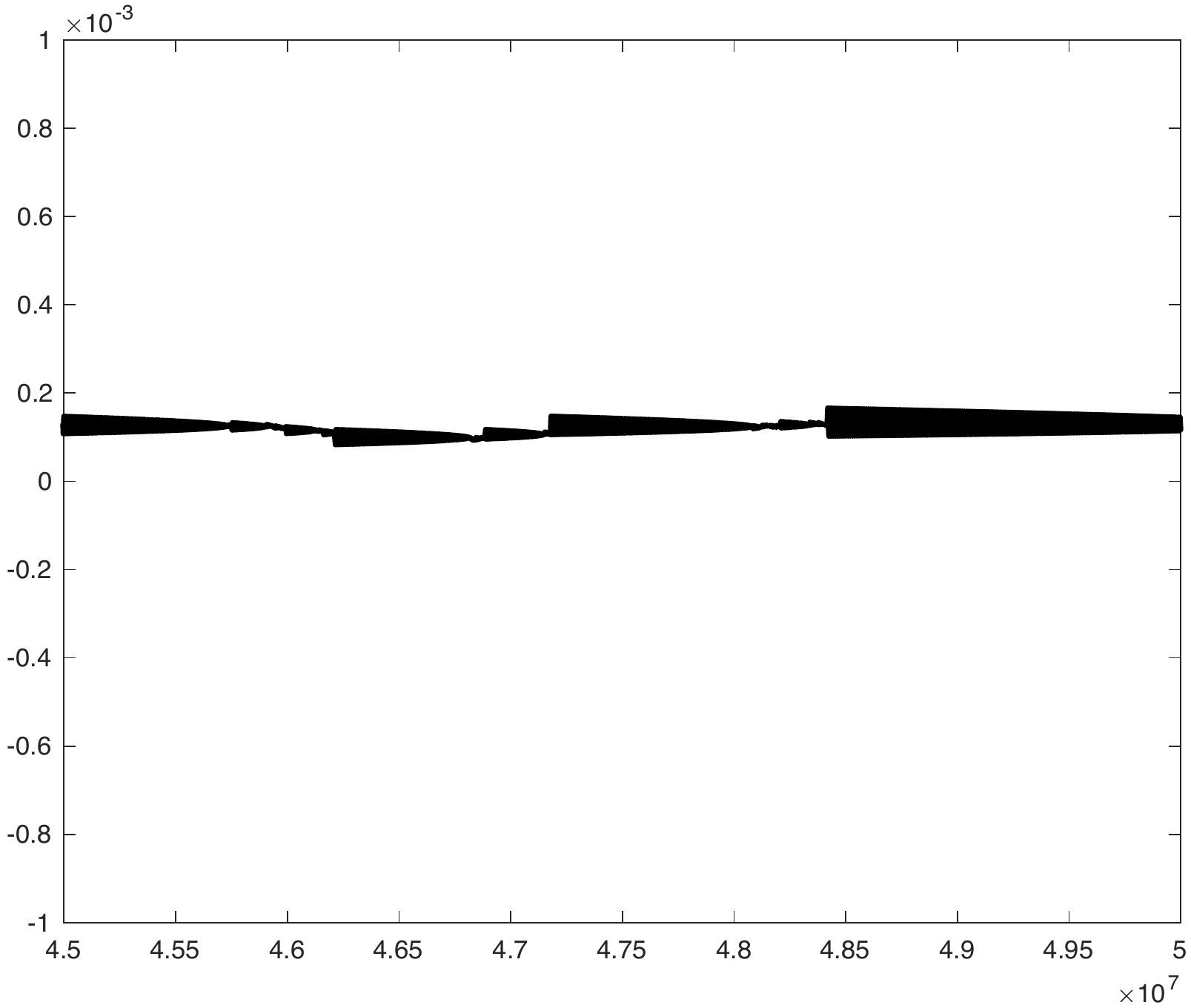}{Same plot as in Fig.~\ref{beta_35}, for the case
$\beta=0.52$.}

\figsimul{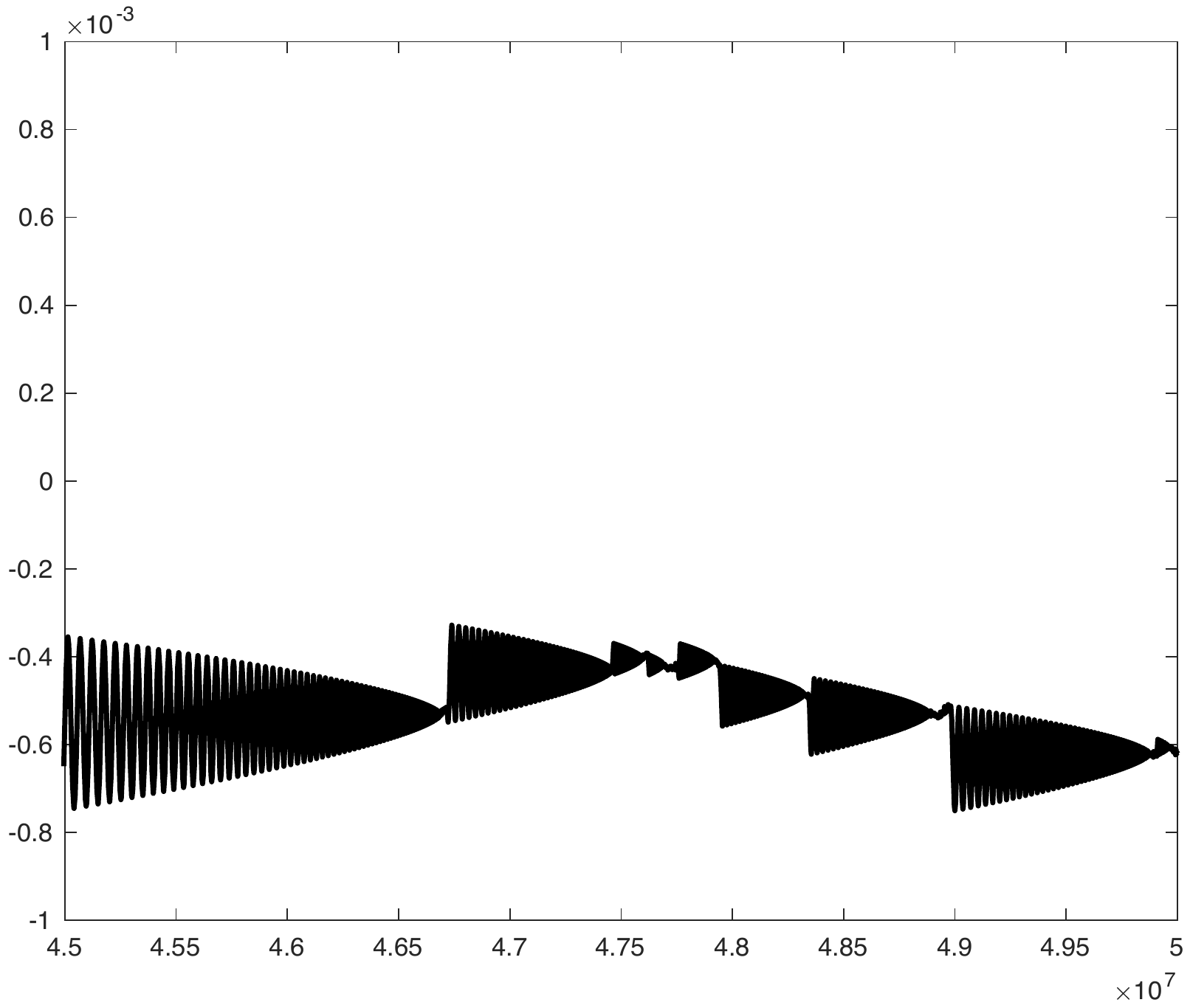}{Same plot as in Fig.~\ref{beta_35}, for the case
$\beta=0.65$.}

We also point the reader to Fig.~\ref{beta_98_scale_1_llong}, which
shows the erratic behavior of Re($\ba_n f$) for large $\beta$. Also,
compare this figure to Fig.~\ref{farey_scale_1_llong}, which displays
the same plot as in Fig.~\ref{beta_98_scale_1_llong} but for the Farey
map $T_F$ as in (\ref{farey-r}). The Farey map is akin to $T_\beta$
with $\beta = 1$. In fact, as can be calculated easily, the partition
$\{ L_k \}$ for $T_F$ is given by $L_0 = [0, \ln 2)$ and $L_k =
[ \ln(k+1) , \ln(k+2) )$ for $k \ge 1$. Thus $\leb( L_k ) \sim k^{-1}$.

\figsimul{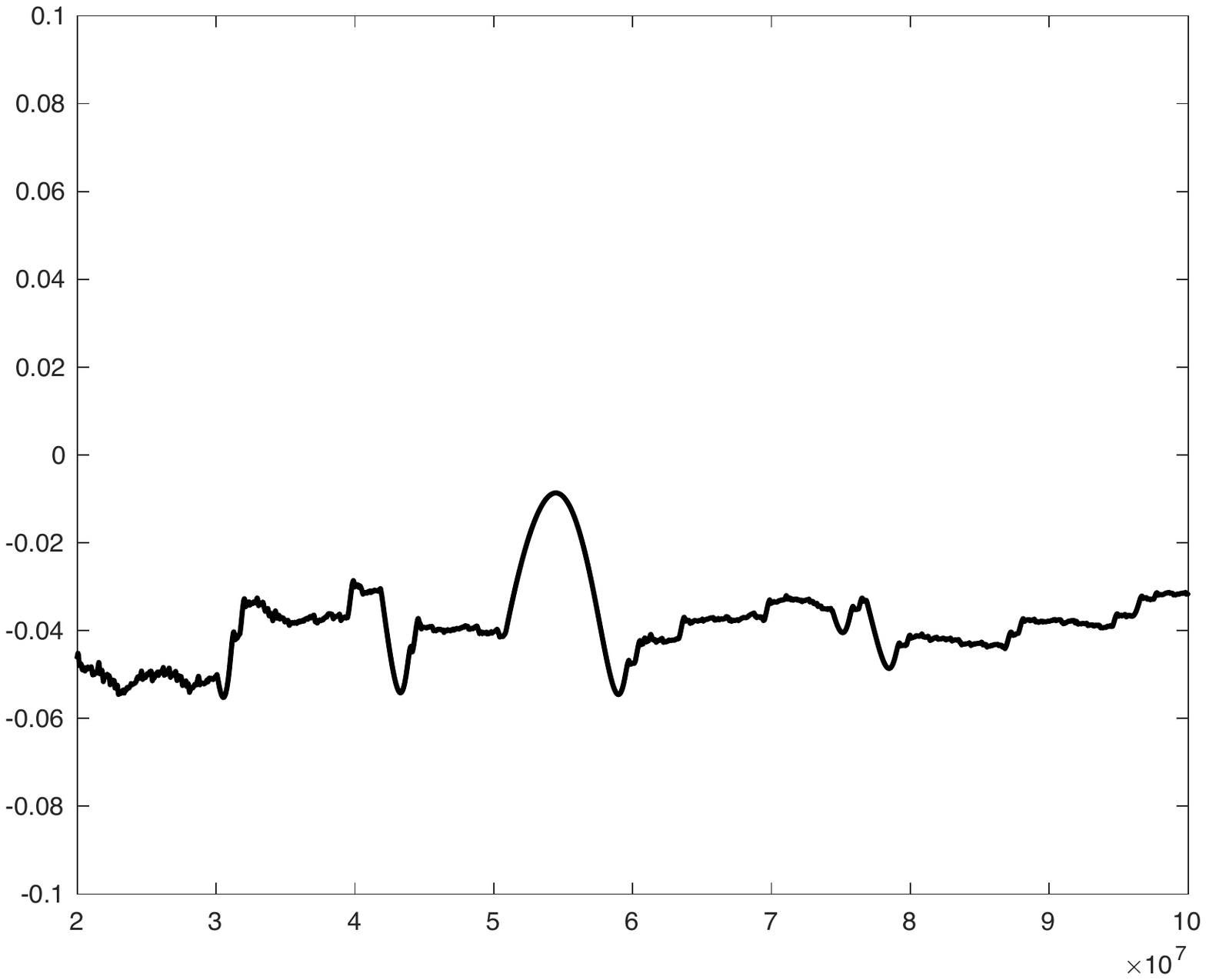}{Plot of $\ba_n g(x_0)$, for the same
$g$ and $x_0$ as in Fig.~\ref{beta_35}, relative to $T_\beta$ with
$\beta=0.98$. Here $0.2 \cdot 10^8 \le n \le 10^8$ and the vertical scale
is in absolute units.}

\figsimul{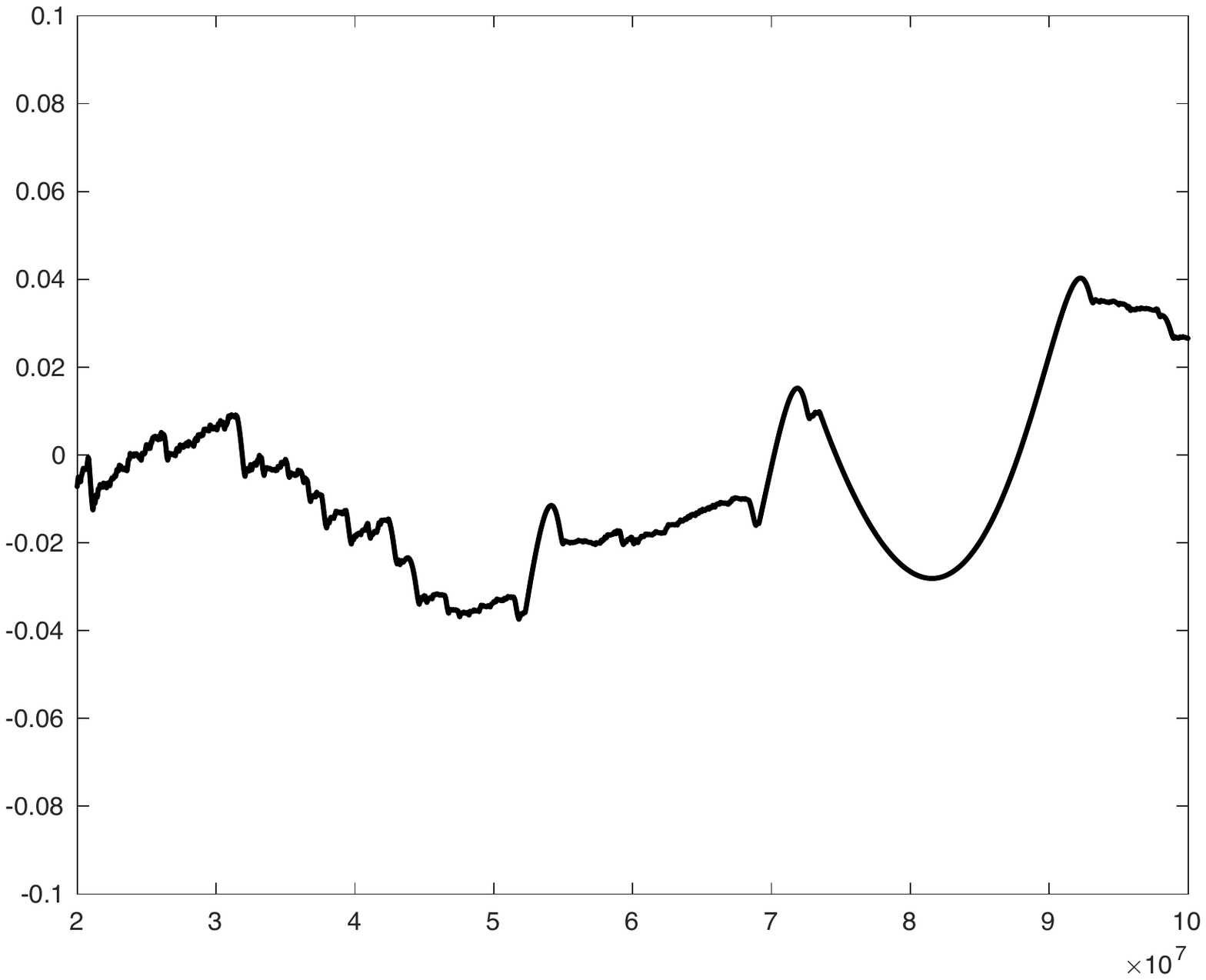}{Plot of $\ba_n g(x_0)$, for the same
$g$ and $x_0$ as in Fig.~\ref{beta_35}, relative to the Farey map $T_F$.
The horizontal range and vertical scale are the same as in
Fig.~\ref{beta_98_scale_1_llong}.}

We can produce more counterexamples to a general Birkhoff Theorem
for $L^\infty$ observables than mentioned above. The ones that we
present momentarily are interesting not only because they do not
follow directly from the results of \cite{sy}, but also because the
Birkhoff sums that we write are representations of \emph{L\'evy walks}.
L\'evy walks are well-studied stochastic processes, often used in nonlinear
and statistical physics as models for anomalous diffusion and transport
\cite{zdk}. In fact, we will use this representation to derive a very fine
limit theorem for our \ob s, thus adding to the connections
between the field of anomalous stochastic processes and infinite
\erg\ theory; cf.\ \cite{k} and references therein.

Our \dsy\ is a Kakutani tower, for which we employ the notation
$(Y, \nu, \tow)$ of Appendix \ref{app-kakutani}. We start by defining
the base map $S: \Sigma \into \Sigma$, where $\Sigma := [0,1) \times
[0,1)$ and $x = (x_1,x_2)$ is a generic element of $\Sigma$.
Let $\mathcal{B} = \{ B_i \}_{i \in \Z^+}$ be a partition of $[0,1)$ made up
of right-open intervals, which are ordered from left to right. Assume also
that $\leb( B_i ) \sim c\, i^{-\beta-1}$, for some $c>0$ and $\beta \in (0,1)$.
Let us define
$S_\mathcal{B}: [0,1) \into [0,1)$ to be the full-branched, piecewise-linear
and increasing Markov map relative to $\mathcal{B}$. In other words,
$S_\mathcal{B} |_{B_i}$ maps $B_i$ onto $[0,1)$ with  derivative
$1/ \leb( B_i )$.
It is clear that $S_\mathcal{B}$ preserves the Lebesgue \me\ $\leb$ and
that the partitions $\mathcal{B}, S_\mathcal{B}^{-1} \mathcal{B}, \ldots,
S_\mathcal{B}^{-n} \mathcal{B}, \ldots$ are independent w.r.t.\ $\leb$.
Then let $\mathcal{C} = \{ C_j \}_{j \in \mathbb{J}}$ be another partition
of $[0,1)$ given by right-open intervals. Here $\mathbb{J}$ can be either
$\{ 1, 2, \ldots, N \}$, for some positive integer $N$, or $\Z^+$. Again
let us assume that the intervals $C_j$ are ordered from left to right.
In analogy with the previous case, we denote $S_\mathcal{C}$ the
full-branched, piecewise-linear and increasing Markov map of $[0,1)$
relative to $\mathcal{C}$. This map has the same properties as
$S_\mathcal{B}$.
Define $S := S_\mathcal{B} \times S_\mathcal{C}$, i.e., $S(x_1,x_2) :=
(S_\mathcal{B}  (x_1), S_\mathcal{C} (x_2))$. Thus $S$ is a
two-dimensional uniformly expanding map which preserves the
Lebesgue \me\ of $\Sigma$; in accordance with the notation of
Appendix \ref{app-kakutani}, this \me\ will be called $\rho$. Also,
the partition $\mathcal{B} \otimes \mathcal{C} := \{ B_i \times C_j \}$ of
$\Sigma$ has the property that all its back-iterates
$S^{-n} ( \mathcal{B} \otimes \mathcal{C} )$ are mutually independent.

For $i \in \Z^+$ denote $A_i := B_i \times [0,1)$. The height \fn\
$\varphi: \Sigma \into \N$ is defined by the identities
\begin{equation}
  \varphi |_{A_i} \equiv i-1.
\end{equation}
Thus $\rho( \{ \varphi \ge k \} ) = \sum_{i > k} \rho (A_i) = \sum_{i > k}
\leb(B_i) \sim c \beta k^{-\beta}$. It follows that the invariant \me\ $\nu$
of $\tow$, which is the Lebesgue \me\ on each level of the tower
\begin{equation} \label{lw-20}
  L_k := \rset{x \in \Sigma} {\varphi(x) \ge k} \times \{k\} = \bigcup_{i \ge k+1}
  A_i \times \{k\} ,
\end{equation}
cf.\ (\ref{def-kt-3}), is infinite.

Lastly, we introduce the \ob\ $f: Y \into \C$. Let
$\{ \gamma_j \}_{j \in \mathbb{J}}$ be a set of complex numbers with
$| \gamma_j | = 1$ and define $f$ so that
\begin{equation} \label{my-ob}
  f |_{[0,1) \times C_j} \equiv \gamma_j.
\end{equation}
The easiest example of such an \ob\ is when $\mathbb{J} = \{ 1, 2 \}$ and
$\gamma_1 = -1$, $\gamma_2 = 1$.

\begin{proposition} \label{prop-lw}
  For the \dsy\ $(Y, \nu, \tow)$ introduced above and the \fn\ $f$ defined
  by (\ref{my-ob}), let us interpret the Birkhoff sum $\bs_n f$ as a random 
  variable for the probability \me\ $\nu_0 := \nu( \,\cdot\, | L_0)$, where 
  this means that
  $\bs_n f(y)$ depends on the initial condition $y = (x,0)$, where $x$
  is chosen randomly in $\Sigma$ according to the Lebesgue \me.
  Then the process
  \begin{displaymath}
    \left( L_n(t) \right)_{t \in \R_0^+} := \left( \frac { \bs_{\lfloor nt \rfloor} f }
    n \right)_{t \in \R_0^+}
  \end{displaymath}
  converges in distribution, w.r.t.\ the topology of the uniform
  convergence on all intervals $[0,T]$, to a continuous $\C$-valued
  process $( L(t) )_{t \in \R_0^+}$.
  If not all $\gamma_j$ are equal (assuming that $\leb(C_j) > 0$ for
  all $j \in \mathbb{J}$), then, for every $t\ge 0$, $L(t)$ is almost surely
  non-constant. In particular
  \begin{displaymath}
    \ba_n f = \frac{ \bs_n f } n
  \end{displaymath}
  converges in distribution to a non-constant random variable.
\end{proposition}

\proof We claim that
\begin{equation}
  \bs_n (f \circ \tow) = \sum_{k=1}^n f \circ \tow^k
\end{equation}
is the \emph{L\'evy walk} on $\C$ thus defined: A walker stands at
the origin of $\C$ when she reads the value of a random integer $I_1$
and a random complex number $\Gamma_1$, with the following
probabilities:
\begin{equation}
  \forall i \in \Z^+, \, j \in \mathbb{J}, \quad \mathrm{Prob}
  \{ (I_1, \Gamma_1) = (i, \gamma_j) \} = \leb(B_i) \,
  \leb(C_j).
\end{equation}
All other values of $(I_1, \Gamma_1)$ occur with probability zero.
Observe that $I_1$ and $\Gamma_1$ are independent by definition;
remember also that $|\Gamma_1| = 1$. The walker then takes $I_1$
unit steps in the direction $\Gamma_1$ \emph{one step at a time} ---
which is why we speak of `walk' instead of `jump'. At this point the
walker reads the value of another random pair $(I_2, \Gamma_2)$,
with the same probabilities as the previous pair and independent of it.
This will determine, in the same way as before, the motion of the
walker during the next $I_2$ time units. And so on.

In other words, we have described a persistent random walk on $\C$,
equivalently, a random walk with an internal state \cite{cgls}, with
long-tailed inertial segments, since $\mathrm{Prob} (I_1 = i) \sim c\,
i^{-\beta-1}$, with $0 < \beta < 1$. (The simple case
$\mathbb{J} = \{ 1, 2 \}$, $\gamma_1 = -1$, $\gamma_2 = 1$
corresponds to a \emph{simple symmetric L\'evy walk} on the real line.)

The claim is not hard to show. Let us for the moment suppose
that the reference, or \emph{initial}, \me\ is not $\nu_0$ but
$\nu( \,\cdot\, | B_{i_0} \times C_{j_0} \times \{0\} )$, for some choice of
$i_0 \in \Z^+$ and $j_0 \in \mathbb{J}$. Using the fact that the
base map $S$ sends each $B_i \times C_j$ affinely onto $\Sigma$
and recalling the definition (\ref{def-kt-4}) of the tower map $\tow$, we
see that the push-forward of the initial \me\ is given by
\begin{equation} \label{lw-100}
  \tow_* \, \nu( \,\cdot\, | B_{i_0} \times C_{j_0} \times \{0\} ) =
  \sum_{i_1, j_1} \leb(B_{i_1}) \,
  \leb(C_{j_1}) \, \nu( \,\cdot\, | B_{i_1} \times C_{j_1} \times
  \{ i_1 - 1 \} ) ,
\end{equation}
where the sum is over $(i_1, j_1) \in \Z^+ \times \mathbb{J}$. If we
we fix one such pair $(i_1, j_1)$ and condition the above to
$B_{i_1} \times C_{j_1}
\times \{ i_1 - 1 \}$ --- more precisely, if we condition the initial \me\ to
the event $\{ \tow(x,0) \in B_{i_1} \times C_{j_1} \times \{ i_1 - 1 \} \} =
\{ S(x) \in B_{i_1} \times C_{j_1} \}$ --- we can push-forward
the resulting \me\ down the levels of the tower. In formula, for all
$1 \le k \le i_1$,
\begin{equation} \label{lw-110}
\begin{split}
  & \tow_*^k \, \nu \! \left( \,\cdot\, | ( B_{i_0} \times C_{j_0} \times \{0\} )
  \cap \tow^{-1} ( B_{i_1} \times C_{j_1} \times \{ i_1 - 1 \} ) \right) \\
  &\quad = \nu \! \left( \,\cdot\, | B_{i_1} \times C_{j_1} \times \{ i_1 - k \}
  \right).
\end{split}
\end{equation}
Therefore, for any $y$ as specified by the conditioning in the above l.h.s.,
$f ( \tow^k (y) ) = \gamma_{j_1}$, whence $\bs_k (f \circ \tow) (y) =
k \gamma_{j_1}$. At time $k = i_1$ the r.h.s.\ of (\ref{lw-110}) is the
Lebesgue \me\ on $B_{i_1} \times C_{j_1} \times \{ 0 \}$, that is, it has
the same form as the initial \me. In other words the process has
\emph{renewed}, losing all memory of the initial \me.

In more detail, this implies that if we fix $q\in \Z^+$, $(i_1, j_1), \ldots,
(i_q, j_q) \in \Z^+ \times \mathbb{J}$ and consider all the initial conditions
$y$ such that the first excursion down the tower starts in
$B_{i_1} \times C_{j_1} \times \{ i_1 - 1 \}$, the second excursion
starts in $B_{i_2} \times C_{j_2} \times \{ i_2 - 1 \}$ and so on up to the
$q^\mathrm{th}$ excursion, then for all $n := i_1 + \cdots + i_{q-1} + k$,
with $1 \le k \le i_q$, we have
\begin{equation}
  \bs_n (f \circ \tow) (y) = i_1 \gamma_{j_1} + \cdots + i_{q-1}
  \gamma_{j_{q-1}} + k \gamma_{j_q}.
\end{equation}
Conditioning to the set of all such $y$ and recalling that the pairs
$(i_1, j_1), \ldots, (i_q, j_q)$ are i.i.d.\ for the initial \me\
$\nu( \,\cdot\, | B_{i_0} \times C_{j_0} \times \{0\} )$ proves our claim, at
least for such choice of the initial \me.

Extending the proof to the case where the initial \me\ is $\nu_0$, as
defined in the statement of the proposition, is immediate. Indeed, the
arguments described above depend in no way on $i_0, j_0$
and $\nu_0$ is a convex linear combination of the probability \me s
$\nu( \,\cdot\, | B_{i_0} \times C_{j_0} \times \{0\} )$, for
$(i_0, j_0) \in \Z^+ \times \mathbb{J}$.

Having established the claim, the assertion of Proposition \ref{prop-lw}
follows from Corollary 4.14 of \cite{mssz} and the fact that
$\bs_n f = f + \bs_{n-1}(f \circ \tow)$. The process $(L(t))$ is a combination
of certain L\'evy processes whose marginals at any fixed time $t$ are
non-constant with probability 1 \cite[Eqs.\ (4.13), (3.10)]{mssz}. (In truth,
the results of \cite{mssz} are stated for the case where $(L_n(t))$ is a
continuous-time process, that is, the walker moves continuously with
unit speed from one ``renewal point'' to the next. Extending such results
to our case is trivial.)
\qed

\begin{remark}
  One might wonder why, in Proposition \ref{prop-lw}, the scaling rate
  of $\bs f_n$ (that is, $n$, a.k.a.\ \emph{ballistic scaling}) does not
  depend on $\beta$, the exponent of the tail of the distribution of the
  inertial segments, when $0 < \beta < 1$. This is a fact about
  L\'evy walks, a rigorous proof of which can be found in \cite{mssz}.
  Here we give a simple, heuristic, explanation. If $\{ \mathcal{X}_n
  \}_{n \in \N}$ denotes the L\'evy walk in the proof of Proposition
  \ref{prop-lw}, let $\{ \mathcal{Y}_k \}_{k \in \N}$ denote its associated
  \emph{L\'evy flight}, defined by $\mathcal{Y}_0 \equiv 0$ and
  $\mathcal{Y}_k := \sum_{q=1}^k I_q \Gamma_q$. (Recall that
  $\{ I_q \}_{q \in \Z^+}$ and $\{ \Gamma_q \}_{q \in \Z^+}$ are two
  independent i.i.d.\ processes such that $I_1$ takes values in
  $\Z^+$ and is in the normal basin of attraction of a skewed
  $\beta$-stable distribution, and $\Gamma_1$ takes values in
  $\mathbb{S}^1 \subset \C$.) In other words, $\{ \mathcal{Y}_k \}$
  is the L\'evy walk $\{ \mathcal{X}_n \}$ seen at its renewal times.
  Furthermore, $\{ \mathcal{X}_n \}$ is a unit-speed interpolation of
  $\{ \mathcal{Y}_k \}$. Now, let $\tau_k := \sum_{q=1}^k I_q$ denote 
  the sequence of renewal times, with $\tau_0 \equiv 0$. By the 
  hypothesis on $I_1$, $\tau_k \approx k^{1/\beta}$, as $k \to \infty$ 
  \cite{il}. The same hypothesis shows that the whole process 
  $\{ k^{-1/\beta} \, \mathcal{Y}_{\lfloor ks \rfloor} \}_{s \in \R_0^+}$
  converges to a L\'evy process $\{ \overline{\mathcal{Y}}(s)
  \}_{s \in \R_0^+}$, in a sense that we do not specify here. Now,
  denote by $n \mapsto  K_n$ the generalized inverse of
  $k \mapsto \tau_k$, i.e., the non-decreasing \fn\ $\N \into \N$ such
  that $\tau_{K_n} \le n < \tau_{K_n + 1}$. Clearly $K_n
  \approx n^\beta$. By construction, $\mathcal{X}_n$ always
  lies between $\mathcal{Y}_{K_n}$ and $\mathcal{Y}_{K_n + 1}$.
  These two processes are not the same --- they are sometimes
  called the \emph{lagging} and \emph{leading} walks of
  $\mathcal{X}_n$, respectively --- but it is easy to show that they
  scale in the same way. For the purposes of this explanation we can
  approximate $\mathcal{X}_n$ with $\mathcal{Y}_{K_n}$, therefore,
  in non-rigorous notation, we can write that, for $n \to \infty$,
  \begin{equation}
    \mathcal{X}_{\lfloor nt \rfloor} \approx
    \mathcal{Y}_{K_{\lfloor nt \rfloor}} \approx
    \mathcal{Y}_{\lfloor n^\beta t^\beta \rfloor} \approx n \,
    \overline{\mathcal{Y}}(t^\beta) .
  \end{equation}
  We end this remark by observing that this explanation only holds
  for $\beta \in (0,1)$. In all other cases, since the first moment of
  $I_1$ is finite (or barely infinite), the scaling of $\mathcal{Y}_k$ is
  generally different from that of $\tau_k$.
\end{remark}

The example that we have presented in Proposition \ref{prop-lw}
is in the same spirit as the occupation times of dynamically
separated sets. In fact, if we denote
\begin{equation}
  R_j := \bigcup_{k \ge 1} \bigcup_{i \ge k+1} B_i \times C_j
  \times \{ k \} ,
\end{equation}
cf.~(\ref{lw-20}), we realize that, for $j \in \mathbb{J}$, the
infinite-\me\ sets $R_j$ are dynamically separated by the juncture
$L_0$. Observe that, in the case where the $\gamma_j$ are all
different, $R_j = (Y \setminus L_0) \cap \{ f = \gamma_j \}$.
In any case, $f$ can be expressed as
\begin{equation}
  f = \sum_{j \in \mathbb{J}} \gamma_j \, 1_{R_j} +
  \sum_{j \in \mathbb{J}} \gamma_j \,1_{[0,1) \times C_j \times
  \{0\} } ,
\end{equation}
where the second sum above amounts to an integrable \fn, cf.\
(\ref{ob-rays}).

Nevertheless, the statistical properties of $\ba_n f$ cannot be
derived from the main theorem of \cite{sy}, and not only because
here we have infinitely many rays. The most important difference
is that the assumption on the asymptotic entrance densities
\cite[Ass.~2.3]{sy} is not satisfied. This follows from the triviality
of the dynamics on the non-zero levels of the tower.

Moreover, our \sy\ can be generalized to the case of uncountably
many rays. It suffices to replace the base map with $S :=
S_\mathcal{B} \times \sigma$, where $S_\mathcal{B}$ is the map
defined earlier and $\sigma$ is the left shift on the space
$( [0,1)^\N, \leb^\N )$ of sequences of i.i.d.\ numbers uniformly
distributed in $[0,1)$. If we define
\begin{equation}
  f(y) = f(x_1, ( \theta_q )_{q \in \N}, k) := e^{2\pi i \theta_0}
\end{equation}
and use the reference \me\ $\nu_0 = \nu( \,\cdot\, | L_0)$
(here $\nu_0$ is isomorphic to $\leb \times \leb^\N$ on $\Sigma :=
[0,1) \times [0,1)^\N$), we see that during the $q^\mathrm{th}$
excursion in the tower the value of $f$ is $e^{2\pi i \theta_q}$ and it is
independent of the values taken during the previous excursions.
Therefore the process $\bs_n f$ is a radially symmetric L\'evy walk on
$\C$. The assertions of Proposition \ref{prop-lw} still hold. Finally,
writing $\N = \{0\} \times \Z^+$, it is clear that the sets
\begin{equation}
  R_{\theta_0} := (Y \setminus L_0) \cap \{ f = e^{2\pi i \theta_0} \} =
  \bigcup_{k \ge 1} \bigcup_{i \ge k+1} B_i \times \left( \{ \theta_0 \}
  \times [0,1)^{\Z^+} \right) \times \{ k \} ,
\end{equation}
for $\theta_0 \in [0,1)$, are dynamically separated rays.

\appendix

\section{Kakutani towers}
\label{app-kakutani}

Let us briefly recall the definition and basic properties of a Kakutani
tower. For more details we refer to \cite[\S1.5]{a}. In this appendix we
restore the indication of the $\sigma$-algebra in the notation. So,
let $(\Sigma, \mathscr{B}, \rho, S)$ be a conservative, non-singular
\dsy\ on a $\sigma$-finite \me\ space, and suppose that
$\varphi: \Sigma \into \N$ is a measurable \fn. Then the tower over
$S$ with \emph{height \fn} $\varphi$ is the \dsy\
$(Y, \mathscr{C}, \nu, \tow)$ defined as follows:
\begin{align}
  & Y := \lset{ (x, k) \in \Sigma \times \N} { 0 \le k \le \varphi(x) } ;
  \label{def-kt-1} \\
  & \mathscr{C} := \sigma \left( \lset{ A \times \{k\} } {k \in \N, \, A
  \subseteq \{\varphi \ge k\} ,\, A \in \mathscr{B} } \right) ;
  \label{def-kt-2} \\
  & \nu(A \times \{k\}) := \rho(A) ; \label{def-kt-3} \\[3pt]
  & \tow(x, k):= \left \{ \begin{array}{ll}
      (x, k-1), & \mbox{if } k \ge 1; \\[2pt]
      (S(x), \varphi(S(x))), & \mbox{if } k=0 .
    \end{array} \right.
  \label{def-kt-4}
\end{align}
(In (\ref{def-kt-2}) the notation $\sigma( \,\cdot\, )$ denotes the
$\sigma$-algebra generated by the sets between parentheses.)
The \emph{tower map} $\tow$ is conservative and non-singular, and if
$\rho \circ S^{-1} = \rho$ then $\nu \circ \tow^{-1} = \nu$. Furthermore,
if $S$ is ergodic, then $\tow$ is ergodic.

Now, suppose that $(\ps, \sa, \mu, T)$ is an invertible
\me-preserving \dsy, and let $\Sigma \in \sa$ be a sweep-out set with
$\mu(\Sigma) > 0$. (Note that this implies that $T$ is conservative, by
Maharam's Recurrence Theorem, see \cite[Thm.~1.1.7]{a}.) Denote the
induced map of $T$ on $\Sigma$ by $T_\Sigma: \Sigma \into \Sigma$,
and by $\sa_\Sigma$ and $\mu_\Sigma$, respectively, the restrictions of
$\sa$ and $\mu$ to $\Sigma$. Also set
\begin{equation}
  \varphi (x) := \min \rset{ n \geq0} {T^{-n-1}(x) \in \Sigma}.
\end{equation}
In other words, $\varphi: \Sigma \into \N$ is the first-return \fn\ of $T^{-1}$ to
the set $\Sigma$, minus one unit.

\begin{proposition} \label{prop-inv}
  The tower constructed over the \dsy\ $(\Sigma, \sa_\Sigma, \mu_\Sigma,
  T_\Sigma)$ w.r.t.\ the height \fn\ $\varphi$ is measure-theoretically
  isomorphic to $(\ps, \sa, \mu, T)$.
\end{proposition}

\proof In the following, we shall always restrict ourselves to the
full-measure set of points in $X$ for which $T^n$ is invertible for
every $n \ge1$.

Since $\Sigma$ is a sweep-out set, for almost every $x \in \ps$
there exists a smallest $n \ge 0$ such that $z := T^n(x) \in \Sigma$.
We define the map $\Phi: \ps \into Y$, where $Y$ is the reference space
of the tower as described above, by setting
\begin{equation}
  \Phi(x) := (z, n).
\end{equation}
This map is well-defined because, by construction, the first-return time
of $z$ to $\Sigma$, w.r.t.\ $T^{-1}$ must be strictly larger than $n$,
implying that $\varphi(z) \ge n$. Clearly $\Phi$ is injective, since,
for $x_1 \ne x_2 \in X$, either these two points have different landing
points $z_1, z_2$ in $\Sigma$, or $z_1 = z_2$. In the latter case,
however, the $T$-trajectories of $x_1, x_2$ cannot get to $z_1 = z_2$
after the same number $n$ of iterations, otherwise the map $T^n$
would not be invertible. Moreover, for a.e.\ $(z, n)\in Y$, if we set
$x:=T^{-n}(z)$, then $\Phi(x) = (z, n)$, showing that $\Phi$ is also
surjective.

It remains to demonstrate that $T = \Phi^{-1} \circ \tow \circ \Phi$.
Let us consider two cases. First, suppose that $x\in X \setminus
\Sigma$. Then we find our $z=T^n(x)\in \Sigma$ with $n \ge1$. So,
\begin{equation}
  \Phi^{-1} \circ \tow (\Phi(x)) = \Phi^{-1}(\tow(z, n)) = \Phi^{-1}(z, n-1)
  = T(x).
\end{equation}
The second case is where $x\in\Sigma$. Here we have that $\Phi(x)
= (z, 0)$, with $z=x$. Thus,
\begin{equation}
  \Phi^{-1} \circ \tow (\Phi(x)) = \Phi^{-1}(\tow(z, 0)) = \Phi^{-1} (
  T_\Sigma(x) , \varphi ( T_\Sigma(x) )) = T(x).
\end{equation}
Here the final equality holds because $\varphi ( T_\Sigma(x) )$ is
equal to $\rho(x) -1$, where $\rho$ denotes the return-time function
to $\Sigma$ with respect to $T$, that is $T_\Sigma(x) := T^{\rho(x)}(x)$.
\qed

Consider again the general set-up from Section \ref{sec:mainthm}. So
$(\ps, \sa, \mu, T)$ is a conservative, ergodic measure-preserving
dynamical system on a $\sigma$-finite, infinite measure space where we
choose a set $L_0 \in \sa$ with $0<\mu(L_0)<\infty$ ($L_0$ is then a
sweep-out set). If it happens that this system is isomorphic to a tower
with $L_0$ identified with the base level $\Sigma \times \{0\}$ (as would
be the case for $T$ invertible, as shown above), then each partition element
$L_k$ is identified with the $k^\mathrm{th}$ level of the tower, i.e.,
$\{ \varphi \ge k \} \times \{k\}$. Moreover, since the tower map sends level
$k$ injectively into level $k-1$, we gain in this case that $T$ maps
$L_k$ injectively into $L_{k-1}$.

\end{document}